\documentclass{CustomArticle}

\usepackage[
  style = trad-abbrv,
  citestyle = numeric-comp,
  sorting = nyt
]{biblatex}
\usepackage{import}
\usepackage{hyperref}
\usepackage[capitalize,compress]{cleveref}

\usepackage{CommonMath}
\usepackage{ProjectMath}
\usepackage{SimpleAlgorithm}
\usepackage{TheoremParts}

\crefname{equation}{}{}

\PaperGeneralInfo{
  title = {Global Complexity Analysis of BFGS},
  date = {April 23, 2024},
}
\AddPaperAuthor{
  name = {Anton Rodomanov},
  affiliation = {CISPA Helmholtz Center for Information Security},
  email = {anton.rodomanov@cispa.de}
}

\PaperAbstract{%
  In this paper, we present a global complexity analysis of the classical
  BFGS method with inexact line search, as applied to minimizing a strongly
  convex function with Lipschitz continuous gradient and Hessian.
  We consider a variety of standard line search strategies including the
  backtracking line search based on the Armijo condition, Armijo--Goldstein and
  Wolfe--Powell line searches.
  Our analysis suggests that the convergence of the algorithm proceeds in
  several different stages before the fast superlinear convergence actually
  begins.
  Furthermore, once the initial point is far away from the minimizer, the
  starting moment of superlinear convergence may be quite large.
  We show, however, that this drawback can be easily rectified by using
  a simple restarting procedure.
}

\begin{document}
  \PrintTitleAndAbstract

  \section{Introduction}

\subsection{Motivation}

Quasi-Newton (QN) methods are classical algorithms for Nonlinear Optimization,
which were invented almost 60 years
ago~\cite{Davidon-VariableMetric-59,Fletcher.Powell-RapidlyConvergent-63} and
became extremely popular due to their excellent practical
performance~\cite{Nocedal.Wright-NumericalOptimization-06}.
QN methods can be seen as intermediate methods between the Gradient and Newton
algorithms, and essentially try to approximate the Newton direction by using
only first-order information about the objective function and at a significantly
smaller cost than that of solving the Newton system of linear equations at each
iteration.
Three popular examples of QN methods are BFGS~\cite{%
  Broyden-ConvergenceClass-70-p1,%
  Broyden-ConvergenceClass-70-p2,%
  Fletcher-NewApproach-70,%
  Goldfarb-FamilyVariable-70,%
  Shanno-ConditioningQuasi-70%
}, SR1~\cite{
  Davidon-VariableMetric-59,%
  Broyden-QuasiNewton-67%
} and DFP~\cite{%
  Davidon-VariableMetric-59,%
  Fletcher.Powell-RapidlyConvergent-63%
} algorithms, among which BFGS is usually considered the most stable and
efficient.

Good practical performance of QN methods quickly stimulated a lot of theoretical
research on these algorithms, most of which was carried out in the 1970--80s.
Among the most important results of that time was the proof of local superlinear
convergence and global convergence of QN methods under both exact and inexact
line search~\cite{%
  Powell-OnTheConvergence-71,%
  DennisMore-ACharacterization-74,
  Powell-SomeGlobal-76,
  Byrd.et.al-GlobalConvergence-87%
}.
However, the Complexity Theory approach, which is nowadays prevalent in
Optimization, was not yet popular back then, and all the corresponding results
were only asymptotic without any specific convergence rate estimates.
For a more detailed discussion on this and additional historical remarks, we
refer the interested reader to~\cite{Rodomanov-QuasiNewton-22}.

Recently, there has been a new wave of theoretical interest in QN algorithms.
Motivated by the desire to obtain concrete efficiency estimates for the local
superlinear convergence of QN methods, \cite{Rodomanov-Nesterov.GreedyQuasi-21}
proposed to slightly modify the standard BFGS, SR1 and DFP methods by replacing
the iterate differences participating in Hessian approximation updates with
greedily selected basis vectors.
The resulting \emph{greedy} QN methods were the first QN algorithms equipped
with explicit convergence rate bounds characterizing their superlinear
convergence.
These methods were later extended to randomized
updates~\cite{Lin.et.al-GreedyRandom-21,Lin.et.al-ExplicitConvergence-22}
possessing improved complexity guarantees, to
solving systems of nonlinear
equations~\cite{Ye.et.al-GreedyRandom-21,Liu.et.al-BlockBroyden-24},
saddle-point problems~\cite{Liu.Luo-QuasiNewton-22}, distributed
optimization~\cite{Du.You-DistributedAdaptive-24},
and a variety of other settings~\cite{%
  Jin.et.al-SharpenedQuasi-22,%
  Liu.et.al-SymmetricRank-23,%
  Ji.Dai-GreedyPSB-23%
},
including QN
variants~\cite{%
  Liu.et.al-IncrementalQuasi-24,%
  Lahoti.et.al-SharpenedLazy-24%
} of the incremental Newton
method~\cite{RodomanovKropotov-SuperlinearlyConvergent-16} for minimizing finite
sums.

At the same time, there have been big advancements in understanding local
complexity guarantees for \emph{classical} QN methods.
The first explicit local convergence analysis of standard BFGS and DFP methods
was done in~\cite{Rodomanov.Nesterov-RatesSuperlinear-22} and subsequently
improved in~\cite{Rodomanov.Nesterov-NewResults-21}.
In particular, it was shown that the standard BFGS method, started from the
neighborhood of the minimizer with an appropriately chosen simple initial
Hessian approximation, requires at most $\BigO(n \log \kappa)$ steps before its
superlinear convergence begins, where $n$ is the space dimension and $\kappa$ is
the condition number of the problem.
These results were later extended in~\cite{Ye.et.al-TowardsExplicit-23} to the
SR1 method utilizing the special correction strategy
from~\cite{Rodomanov-Nesterov.GreedyQuasi-21}.
See also~\cite{Jin.Mokhtari-NonAsymptotic-23} for the analysis of QN methods
under the assumption that the initial Hessian approximation is sufficiently
accurate, and~\cite{Lin.et.al-ExplicitSuperlinear-21} for the local analysis of
Broyden's method for solving nonlinear equations.

While the focus of the previously mentioned works was on \emph{local}
convergence, there has also been some recent research on obtaining
explicit \emph{global} efficiency estimates for QN methods, including
accelerated algorithms~\cite{%
  Jiang.et.al-OnlineLearning-23,%
  Jiang.Mokhtari-24-AcceleratedQuasi,%
  Kamzolov.et.al-CubicRegularization-23,%
  Scieur-AdaptiveQuasi-24%
}.
These works develop new methods with interesting complexity guarantees
but the algorithms being studied are different from standard QN methods,
and either do not have superlinear convergence at all, or their superlinear
convergence and the actual practical behavior might be worse than that of the
classical QN methods.

It is therefore interesting to understand which global complexity guarantees we
have for \emph{classical} QN methods, e.g., the standard BFGS algorithm.
Two papers have appeared very recently addressing this question:
\cite{Krutikov.et.al-ConvergenceRate-23} and~\cite{Jin.et.al-NonAsymptotic-24}.
However, they both study QN methods under a rather strong assumption on the
\emph{exact line search}, whose additional complexity still needs to be properly
estimated and under which many QN methods, including BFGS, DFP and SR1,
coincide~\cite{Dixon-QuasiNewton-72-p1,Dixon-QuasiNewton-72-p2}.

In this work, we make a further advancement in this direction by establishing
new complexity results for the standard BFGS method for a variety of
\emph{inexact line search} strategies which are simple to implement and whose
complexity can be properly accounted for.
\subsection{Overview of Main Results and Contributions}

We consider the problem of unconstrained minimization:
\begin{equation}
  \label{eq:Intro-Problem}
  \min_{x \in \RealField^n} f(x),
\end{equation}
where $\Map{f}{\RealField^n}{\RealField}$ is a $\mu$-strongly convex function
with $L$-Lipschitz gradient and $L_2$-Lipschitz Hessian.
We study the standard BFGS method for solving~\eqref{eq:Intro-Problem}:
\begin{equation}
  \label{eq:Intro-BFGS}
  \begin{aligned}
    x_{k + 1} &= x_k - h_k H_k \Gradient f(x_k),
    \\
    H_{k + 1} &= \InvBFGS(H_k, \delta_k, \gamma_k),
  \end{aligned}
\end{equation}
where $\InvBFGS$ is the inverse BFGS update formula (see \cref{eq:Def-InvBFGS}),
$\delta_k \DefinedEqual x_{k + 1} - x_k$,
$\gamma_k \DefinedEqual \Gradient f(x_{k + 1}) - \Gradient f(x_k)$, and
$h_k \geq 0$ is a ceartainly chosen step size at each iteration~$k \geq 0$.

We show that, for a variety of standard line search strategies,
including the backtracking line search based on the Armijo condition,
Armijo--Goldstein and Wolfe--Powell line searches,
method~\eqref{eq:Intro-BFGS} has global linear convergence
with constant $\kappa^2$, where $\kappa \DefinedEqual L / \mu$ is the condition
number.
This global convergence, however, might not begin immediately but after a
certain preliminary phase whose length is at most $\BigO(n \log \kappa)$ for any
``reasonable'' initial inverse Hessian approximation, namely, the one
satisfying $\frac{1}{L} I \preceq H_0 \preceq \frac{1}{\mu} I$, which includes
the following choices commonly used in practice (see, e.g., Section~6.1
in~\cite{Nocedal.Wright-NumericalOptimization-06}):
$H_0 = \frac{\InnerProduct{\gamma_0'}{\delta_0'}}{\Norm{\gamma_0'}^2}$
or $H_0 = \frac{\Norm{\delta_0'}^2}{\InnerProduct{\gamma_0'}{\delta_0'}} I$,
where $\delta_0' \DefinedEqual x_0' - x_0$,
$\gamma_0' \DefinedEqual \Gradient f(x_0') - \Gradient f(x_0)$,
and $x_0'$ is an arbitrary point distinct from~$x_0$.

We then prove that, after at most
$K_0^L = \BigO(\kappa^2 \log[H^2 F_0] + n \log \kappa)$
steps, where $F_0 \DefinedEqual f(x_0) - f^*$ is the initial function residual
and $H$ is the constant of self-concordancy ($H \leq L_2 / \mu^{3 / 2}$),
method~\eqref{eq:Intro-BFGS} arrives into the neighborhood of local superlinear
convergence.
However, in order for the superlinear convergence to actually begin (after which
it takes no more than a very small number, $\log(1/[H^2 \epsilon])$, iterations
to converge to an $\epsilon$-approximate solution in terms of function
residual), the method needs to pass through another auxiliary phase whose length
depends on the initial Hessian approximation $H_0$, as well as the initial
function residual~$F_0$.
In particular, when the starting point~$x_0$ was not initially in the
neighborhood of the local superlinear convergence, this auxiliary phase may last
much longer than $K_1^L = \BigO(n \log \kappa)$, which is the best known so far
starting moment of superlinear convergence of BFGS established
in~\cite{Rodomanov.Nesterov-NewResults-21}.

To mitigate the possibility of slow superlinear convergence when the initial
function residual~$F_0$ is large, we propose a simple restarting procedure
which prescribes restarting the algorithm every $n, 2 n, 4 n, 8 n, \ldots$
iterations with the fixed initial Hessian approximation~$H_0 = \frac{1}{\mu} I$.
Equipped with this procedure, method~\eqref{eq:Intro-BFGS} still enjoys the
same worst-case global complexity guarantee of $K_0^L$ to reach the region of
local superlinear convergence, as well as the best known so far
$\BigO(n \log \kappa)$ complexity to start the fast superlinear convergence.
As a result, in total, the algorithm needs at most the following number of
iterations to obtain an $\epsilon$-approximate solution to
problem~\eqref{eq:Intro-Problem}:
\[
  \BigO\bigl(
    \kappa^2 \log[H^2 F_0] + n \log \kappa + \log(1/[H^2 \epsilon]
  \bigr).
\]

It is worth mentioning that, for method~\eqref{eq:Intro-BFGS}, we cannot
simply use the already existing complexity results on local superlinear
convergence of BFGS such as those
from~\cite{%
  Rodomanov.Nesterov-NewResults-21,%
  Rodomanov.Nesterov-RatesSuperlinear-22%
}.
The main problem is that all such results are valid only for the algorithm
using unit step sizes $h_k \equiv 1$ at each iteration.
Even though one can indeed prove that any of the line search strategies we
consider will indeed eventually switch to the unit step size, one cannot,
unfortunately, easily obtain a good bound on the corresponding moment when the
switch occurs.
An important contribution of this work is the extension of the local analysis
from~\cite{Rodomanov.Nesterov-RatesSuperlinear-22}, which now requires only a
mild upper bound on the average of $\Abs{h_k - 1}$ across iterations instead of
the explicit assumption of unit step sizes.

\subsection{Notation}

Everywhere in this paper, $\VectorSpace{E}$ denotes an arbitrary
finite-dimensional real vector space.
Its dual space, composed of all linear functionals on~$\VectorSpace{E}$,
is denoted by~$\VectorSpace{E}\Dual$.
The value of a linear functional~$s \in \VectorSpace{E}\Dual$,
evaluated at a point~$x \in \VectorSpace{E}$, is denoted by
$\DualPairing{s}{x}$.

Any positive definite linear operator~$C$ from~$\VectorSpace{E}$
to~$\VectorSpace{E}\Dual$
(notation: $C \in \PDLin(\VectorSpace{E}, \VectorSpace{E}\Dual)$) induces the
following Euclidean norm in~$\VectorSpace{E}$:
\begin{equation}
  \label{Def-RelativeNorm}
  \RelativeNorm{u}{C} \DefinedEqual \DualPairing{C u}{u}^{1 / 2},
  \qquad
  \forall u \in \VectorSpace{E}.
\end{equation}
The corresponding conjugate norm is
\begin{equation}
  \label{RelativeDualNorm}
  \RelativeDualNorm{s}{C}
  \DefinedEqual
  \max_u \SetBuilder{\DualPairing{s}{u}}{\RelativeNorm{u}{C} = 1}
  =
  \DualPairing{s}{C^{-1} s}^{1 / 2},
  \qquad
  \forall s \in \VectorSpace{E}\Dual.
\end{equation}
In the case when $C = \Hessian f(x)$ for a twice differentiable function
$\Map{f}{\VectorSpace{E}}{\RealField}$ and a point $x \in \VectorSpace{E}$, we
prefer to use the notation $\RelativeNorm{\cdot}{x}$ and
$\RelativeDualNorm{\cdot}{x}$, provided there is no ambiguity with the
function~$f$.

It will also be convenient to introduce a fixed operator
$B \in \PDLin(\VectorSpace{E}, \VectorSpace{E}\Dual)$ and denote
\begin{equation}
  \label{Def-Norms}
  \Norm{u} \DefinedEqual \RelativeNorm{u}{B},
  \qquad
  \DualNorm{s} \DefinedEqual \RelativeDualNorm{s}{B},
  \qquad
  \forall u \in \VectorSpace{E}, \
  \forall s \in \VectorSpace{E}\Dual.
\end{equation}

For any pair of linear operators
$\Map{H}{\VectorSpace{E}\Dual}{\VectorSpace{E}}$ and
$\Map{A}{\VectorSpace{E}}{\VectorSpace{E}\Dual}$, we define their \emph{trace}
and \emph{determinant products} by
$\DualPairing{H}{A} \DefinedEqual \Trace(H\Adjoint A)$
and $\det(H, A) \DefinedEqual \det(H A)$, respectively.

For more details on dual spaces and other linear algebra notions we have
mentioned above, we refer the reader to Section~2.1
in~\cite{Rodomanov-QuasiNewton-22}.

  \section{BFGS Update and Its Properties}

Given $H \in \PDLin(\VectorSpace{E}\Dual, \VectorSpace{E})$ and
$\delta \in \VectorSpace{E}$, $\gamma \in \VectorSpace{E}\Dual$, such that
$\DualPairing{\gamma}{\delta} > 0$, we define the \emph{BFGS update} of~$H$
w.r.t.\ the pair $(\gamma, \delta)$ in the usual way:
\begin{equation}
  \label{eq:Def-InvBFGS}
  \InvBFGS(H, \delta, \gamma)
  \DefinedEqual
  H
  -
  \frac{
    H \gamma \delta\Adjoint + \delta \gamma\Adjoint H
  }{
    \DualPairing{\gamma}{\delta}
  }
  +
  \biggl(
    \frac{\DualPairing{\gamma}{H \gamma}}{\DualPairing{\gamma}{\delta}}
    +
    1
  \biggr)
  \frac{\delta \delta\Adjoint}{\DualPairing{\gamma}{\delta}}.
\end{equation}
It is well-known that the result of this update is guaranteed to be positive
definite.

This update corresponds to the following primal update:
\begin{equation}
  \label{BFGS-Primal}
  \BFGS(G, \delta, \gamma)
  \DefinedEqual
  G
  -
  \frac{G \delta \delta\Adjoint G}{\DualPairing{G \delta}{\delta}}
  +
  \frac{\gamma \gamma\Adjoint}{\DualPairing{\gamma}{\delta}},
\end{equation}
where
$G \DefinedEqual H^{-1} \in \PDLin(\VectorSpace{E}, \VectorSpace{E}\Dual)$.

In our analysis, we will need to measure the distance between positive definite
linear operators.
For this, let us introduce the following \emph{prox function}
$\Map{d}{\PDLin(\VectorSpace{E}, \VectorSpace{E}\Dual)}{\RealField}$:
\begin{equation}
  \label{ProxFunc}
  d(X) \DefinedEqual -\ln \det(D^{-1}, X),
\end{equation}
where $D \in \PDLin(\VectorSpace{E}, \VectorSpace{E}\Dual)$ is a certain
(arbitrarily fixed) scaling operator.
Note that the function~$d$ is strictly convex and differentiable.
Therefore, we can use it to define the corresponding
\emph{Bregman distance}
$
  \Map{\psi}{
    \PDLin(\VectorSpace{E}, \VectorSpace{E}\Dual)
    \times
    \PDLin(\VectorSpace{E}, \VectorSpace{E}\Dual)
  }{\RealField}
$:
\begin{equation}
  \label{BregmanDist}
  \begin{aligned}
    \psi(X, Y)
    &\DefinedEqual
    d(Y) - d(X) - \DualPairing{\Gradient d(X)}{Y - X}
    \\
    &=
    -\ln \det(D^{-1}, Y)
    + \ln \det(D^{-1}, X)
    -
    \DualPairing{-X^{-1}}{Y - X}
    \\
    &=
    \DualPairing{X^{-1}}{Y - X}
    -
    \ln \det(X^{-1}, Y).
  \end{aligned}
\end{equation}
Note that, by the strict convexity of the prox function~$d$, for any
$X, Y \in \PDLin(\VectorSpace{E}, \VectorSpace{E}\Dual)$, we have
$\psi(X, Y) \geq 0$ with the equality iff $X = Y$.
Note also that the Bregman distance~\eqref{BregmanDist} does not depend on the
choice of the scaling operator~$D$ in \cref{ProxFunc}.
From \cref{BregmanDist}, for any
$X, Y \in \PDLin(\VectorSpace{E}, \VectorSpace{E}\Dual)$,
we also obtain the following representation of $\psi(X, Y)$ via the relative
eigenvalues $\lambda_1 \geq \cdots \geq \lambda_n > 0$ of~$Y$ w.r.t.~$X$:
\begin{equation}
  \label{BregmanDist-ViaEigs}
  \psi(X, Y)
  =
  \sum_{i = 1}^n (\lambda_i - 1 - \ln \lambda_i)
  =
  \sum_{i = 1}^n \omega(\lambda_i - 1),
\end{equation}
where $\Map{\omega}{\OpenOpenInterval{-1}{+\infty}}{\RealField}$ is the
function
\begin{equation}
  \label{Def-omega}
  \omega(t) \DefinedEqual t - \ln(1 + t).
\end{equation}
Note that $\omega$ is a nonnegative strictly convex function.
It is strictly decreasing on $\OpenClosedInterval{-1}{0}$, strictly increasing
on $\ClosedOpenInterval{0}{+\infty}$, and has a minimizer at~$0$.

Let us establish how the Bregman distance
from~$G \in \PDLin(\VectorSpace{E}, \VectorSpace{E}\Dual)$ to a certain fixed
operator~$C \in \PDLin(\VectorSpace{E}, \VectorSpace{E}\Dual)$
changes after one BFGS update to~$G$.
In what follows, for any
$C \in \PDLin(\VectorSpace{E}, \VectorSpace{E}\Dual)$ and any
$\delta \in \VectorSpace{E}$, $\gamma \in \VectorSpace{E}\Dual$, such that
$\DualPairing{\gamma}{\delta} > 0$, we denote
\begin{equation}
  \label{Def-Alpha-Beta}
  \alpha_C(\delta, \gamma)
  \DefinedEqual
  \frac{\RelativeDualNorm{\gamma}{C}^2}{\DualPairing{\gamma}{\delta}},
  \qquad
  \beta_C(\delta, \gamma)
  \DefinedEqual
  \frac{\DualPairing{\gamma}{\delta}}{\RelativeNorm{\delta}{C}^2}.
\end{equation}
Since $\DualPairing{\gamma}{\delta} > 0$, the fractions
$\alpha_C(\delta, \gamma)$ and $\beta_C(\delta, \gamma)$ are both well-defined
and strictly positive.
Moreover, by the Cauchy--Schwarz inequality, we have
\begin{equation}
  \label{Alpha-Beta-Relation}
  \alpha_C(\delta, \gamma) \geq \beta_C(\delta, \gamma).
\end{equation}
Further, for any $\alpha \geq \beta > 0$, it will be convenient to define
\begin{equation}
  \label{Def-phi}
  \phi(\alpha, \beta)
  \DefinedEqual
  \alpha - 1 - \ln \beta \ (\geq 0).
\end{equation}
Finally, for any $C \in \PDLin(\VectorSpace{E}, \VectorSpace{E}\Dual)$ and any
$\delta \in \VectorSpace{E}$, $\gamma \in \VectorSpace{E}\Dual$, such that
$\DualPairing{\gamma}{\delta} > 0$, let us define
\begin{equation}
  \label{Def-phi-Short}
  \phi_C(\delta, \gamma)
  \DefinedEqual
  \phi\bigl( \alpha_C(\delta, \gamma), \beta_C(\delta, \gamma) \bigr).
\end{equation}
Observe that, for each $C \in \PDLin(\VectorSpace{E}, \VectorSpace{E}\Dual)$,
the quantities $\alpha_C(\cdot, \cdot)$, $\beta_C(\cdot, \cdot)$, and hence
$\phi_C(\cdot, \cdot)$, are homogeneous:
for all $t \in \RealField \setminus \Set{0}$, we have
\begin{equation}
  \label{Alpha-Beta-Phi-Homogeneity}
  \alpha_C(t \delta, t \gamma) = \alpha_C(\delta, \gamma),
  \qquad
  \beta_C(t \delta, t \gamma) = \beta_C(\delta, \gamma),
  \qquad
  \phi_C(t \delta, t \gamma) = \phi_C(\delta, \gamma).
\end{equation}

\begin{lemma}
  \label{St:BregmanDist-Upd}
  Let $C, G \in \PDLin(\VectorSpace{E}, \VectorSpace{E}\Dual)$,
  and let $\delta \in \VectorSpace{E}$, $\gamma \in \VectorSpace{E}\Dual$ be
  such that $\DualPairing{\gamma}{\delta} > 0$.
  Then, for $G_+ \DefinedEqual \BFGS(G, \delta, \gamma)$, we have
  \[
    \psi(C, G) - \psi(C, G_+)
    =
    \phi_C(\delta, G \delta) - \phi_C(\delta, \gamma).
  \]
\end{lemma}

\begin{proof}
  According to \cref{BregmanDist}, we have
  \begin{equation}
    \label{Proof:St:BregmanDist-Upd:Diff}
    \begin{aligned}
      \hspace{1em}&\hspace{-1em}
      \psi(C, G) - \psi(C, G_+)
      \\
      &=
      \bigl[ \DualPairing{C^{-1}}{G - C} - \ln \det(C^{-1}, G) \bigr]
      -
      \bigl[ \DualPairing{C^{-1}}{G_+ - C} - \ln \det(C^{-1}, G_+) \bigr]
      \\
      &=
      \DualPairing{C^{-1}}{G - G_+}
      +
      \ln \det(G, G_+^{-1}).
    \end{aligned}
  \end{equation}

  From \cref{BFGS-Primal,Def-Alpha-Beta}, it follows that
  \begin{equation}
    \label{Proof:St:BregmanDist-Upd:DualPair}
    \DualPairing{C^{-1}}{G - G_+}
    =
    \frac{\RelativeDualNorm{G \delta}{C}^2}{\DualPairing{G \delta}{\delta}}
    -
    \frac{\RelativeDualNorm{\gamma}{C}^2}{\DualPairing{\gamma}{\delta}}
    =
    \alpha_C(\delta, G \delta) - \alpha_C(\delta, \gamma).
  \end{equation}
  Also, it is well-known that
  \[
    \det(G^{-1}, G_+)
    =
    \frac{\DualPairing{\gamma}{\delta}}{\DualPairing{G \delta}{\delta}}
  \]
  (for the proof, see, e.g., Lemma~6.2
  in~\cite{Rodomanov.Nesterov-RatesSuperlinear-22}).
  Hence,
  \begin{equation}
    \label{Proof:St:BregmanDist-Upd:LogDet}
    \ln \det(G^{-1}, G_+)
    =
    \ln \frac{\DualPairing{\gamma}{\delta}}{\DualPairing{G \delta}{\delta}}
    =
    \ln \Bigl(
      \frac{\DualPairing{\gamma}{\delta}}{\RelativeNorm{\delta}{C}^2}
      \frac{\RelativeNorm{\delta}{C}^2}{\DualPairing{G \delta}{\delta}}
    \Bigr)
    =
    \ln \frac{\beta_C(\delta, \gamma)}{\beta_C(\delta, G \delta)},
  \end{equation}
  where the last identity follows from \cref{Def-Alpha-Beta}.

  Substituting
  \cref{Proof:St:BregmanDist-Upd:DualPair,Proof:St:BregmanDist-Upd:LogDet} into
  \cref{Proof:St:BregmanDist-Upd:Diff}, we obtain
  \[
    \begin{aligned}
      \hspace{3em}&\hspace{-3em}
      \psi(C, G) - \psi(C, G_+)
      =
      \alpha_C(\delta, G \delta) - \alpha_C(\delta, \gamma)
      +
      \ln \frac{\beta_C(\delta, \gamma)}{\beta_C(G \delta, \gamma)}
      \\
      &=
      \bigl[
        \alpha_C(\delta, G \delta) - 1 - \ln \beta_C(\delta, G \delta)
      \bigr]
      -
      \bigl[
        \alpha_C(\delta, \gamma) - 1 - \ln \beta_C(\delta, \gamma)
      \bigr]
      \\
      &=
      \phi_C(\delta, G \delta) - \phi_C(\delta, \gamma),
    \end{aligned}
  \]
  where the last identity follows from \cref{Def-phi,Def-phi-Short}.
\end{proof}

  \section{BFGS Method}

Consider the unconstrained optimization problem
\begin{equation}
  \label{Prob}
  \min_{x \in \VectorSpace{E}} f(x),
\end{equation}
where $\Map{f}{\VectorSpace{E}}{\RealField}$ is a differentiable strongly
convex function.
Note that, under our assumptions, problem~\eqref{Prob} has a unique solution.
In what follows, we denote this solution and the corresponding optimal value
by~$x_*$ and~$f_*$, respectively.

Consider the following BFGS Method for solving problem~\eqref{Prob}.

\begin{SimpleAlgorithm}[
    title = {BFGS Method},
    label = Alg:BFGS,
    width = 0.75\linewidth
  ]
  \begin{AlgorithmGroup}[Initialization]
    Choose $x_0 \in \VectorSpace{E}$ and
    $H_0 \in \PDLin(\VectorSpace{E}\Dual, \VectorSpace{E})$.
  \end{AlgorithmGroup}
  \AlgorithmGroupSeparator
  \begin{AlgorithmGroup}[Iteration $k \geq 0$]
    \begin{AlgorithmSteps}
      \AlgorithmStep
        Compute $d_k \DefinedEqual H_k \Gradient f(x_k)$.
      \AlgorithmStep
        Choose a step size $h_k > 0$.
      \AlgorithmStep
        Set $x_{k + 1} \DefinedEqual x_k - h_k d_k$.
      \AlgorithmStep
        \label{Alg:BFGS:Def-delta-gamma}
        Compute $\delta_k \DefinedEqual x_{k + 1} - x_k$ and
        $\gamma_k \DefinedEqual \Gradient f(x_{k + 1}) - \Gradient f(x_k)$.
      \AlgorithmStep
        Set $H_{k + 1} = \InvBFGS(H_k, \delta_k, \gamma_k)$.
    \end{AlgorithmSteps}
  \end{AlgorithmGroup}
\end{SimpleAlgorithm}

Later we will specify more precisely the rules for choosing the step size $h_k$
at each iteration $k \geq 0$ of \cref{Alg:BFGS}.
Also, to avoid considering certain degenerate cases all the time, in what
follows, we will assume that, in \cref{Alg:BFGS},
\[
  \DualPairing{\gamma_k}{\delta_k} > 0
\]
for all $k \geq 0$.
Indeed, if $\DualPairing{\gamma_k}{\delta_k} = 0$ for some $k \geq 0$, then, by
the strong convexity of the objective function and the fact that $H_k$ is
positive definite, we necessarily have $x_{k + 1} = x_k = x_*$, which means
that \cref{Alg:BFGS} has identified an optimal solution.

Note that, in \cref{Alg:BFGS}, for each $k \geq 0$, the operator $H_k$ is
positive definite and hence $d_k$ is a descent direction of~$f$ at~$x_k$:
\[
  \DualPairing{\Gradient f(x_k)}{d_k} \geq 0.
\]

  \section{Global Linear Convergence}

We assume that the objective function in problem~\eqref{Prob} is strongly
convex, and its gradient is Lipschitz continuous: there exist constant
$\mu, L > 0$, such that, for all $x, y \in \VectorSpace{E}$, it holds
\begin{gather}
  \label{StrongConv}
  \DualPairing{\Gradient f(x) - \Gradient f(y)}{x - y}
  \geq
  \mu \Norm{x - y}^2,
  \\
  \label{LipGrad}
  \DualNorm{\Gradient f(x) - \Gradient f(y)}
  \leq
  L \Norm{x - y}.
\end{gather}

Since the function~$f$ is convex, condition~\eqref{LipGrad} is, in fact,
equivalent to any of the following inequalities
for all $x, y \in \VectorSpace{E}$
(see, e.g., Theorem~2.1.5 in~\cite{Nesterov-LecturesConvex-18}):
\begin{gather}
  \label{LipGrad-FuncUbd}
  f(y)
  \leq
  f(x) + \DualPairing{\Gradient f(x)}{y - x} + \frac{L}{2} \Norm{y - x}^2,
  \\
  \label{LipGrad-ViaGrad}
  \DualPairing{\Gradient f(x) - \Gradient f(y)}{x - y}
  \geq
  \frac{1}{L} \DualNorm{\Gradient f(x) - \Gradient f(y)}^2.
\end{gather}
It is also well-known that assumption~\eqref{StrongConv} implies that, for all
$x, y \in \VectorSpace{E}$, we have
\begin{gather}
  \label{StrongConv-FuncLbd}
  f(y)
  \geq
  f(x)
  +
  \DualPairing{\Gradient f(x)}{y - x}
  +
  \frac{\mu}{2} \Norm{y - x}^2,
  \\
  \label{StrongConv-FuncUbd}
  f(y)
  \leq
  f(x)
  +
  \DualPairing{\Gradient f(x)}{y - x}
  +
  \frac{1}{2 \mu} \DualNorm{\Gradient f(y) - \Gradient f(x)}^2.
\end{gather}

An important characteristic of our objective function is its
\emph{condition number}:
\begin{equation}
  \label{ConditionNumber}
  \kappa \DefinedEqual \frac{L}{\mu} \ (\geq 1).
\end{equation}

From now on, we make the following main assumption about the step sizes in
\cref{Alg:BFGS}.

\begin{assumption}
  \label{As-StepSize-Global}
  At each iteration $k \geq 0$ of \cref{Alg:BFGS}, the step size $h_k$ is chosen
  in such a way so that
  \begin{equation}
    \label{FuncProg}
    f(x_k) - f(x_{k + 1})
    \geq
    \eta \frac{\DualPairing{\Gradient f(x_k)}{d_k}^2}{L \Norm{d_k}^2},
  \end{equation}
  where $\eta > 0$ is a certain (absolute) constant.
\end{assumption}

To ensure that \cref{As-StepSize-Global} holds, we can use many line search
strategies.

\begin{example}[``Constant'' step size]
  Suppose the Lipschitz constant~$L$ is known.
  Then, we can choose the step size $h_k$ which minimizes the upper bound on
  $f(x_k - h d_k)$ resulting by using \cref{LipGrad-FuncUbd}:
  \begin{equation}
    \label{LS-LipUbd}
    f(x_k - h d_k)
    \leq
    f(x_k)
    -
    h \DualPairing{\Gradient f(x_k)}{d_k}
    +
    \frac{L}{2} h^2 \Norm{d_k}^2.
  \end{equation}
  This leads us to the step size
  \[
    h_k
    \DefinedEqual
    \frac{\DualPairing{\Gradient f(x_k)}{d_k}}{L \Norm{d_k}^2}.
  \]
  Substituting this expression into \cref{LS-LipUbd}, we obtain \cref{FuncProg}
  with
  \[
    \eta \DefinedEqual \frac{1}{2}.
  \]
\end{example}

\begin{example}[Backtracking line search]\label{St:Backtrack}
  In this strategy, given some
  \begin{equation}
    \label{Backtrack:StartL-Ubd}
    L_{k, 0} \leq L,
  \end{equation}
  we find the smallest index $i_k \geq 0$ over $i = 0, 1, 2, \dots$ so that for
  \begin{equation}
    \label{Backtrack:Quants}
    L_{k, i}
    \DefinedEqual
    2^i L_{k, 0},
    \qquad
    h_{k, i}
    \DefinedEqual
    \frac{\DualPairing{\Gradient f(x_k)}{d_k}}{L_{k, i} \Norm{d_k}^2},
    \qquad
    x_{k + 1, i}
    \DefinedEqual
    x_k - h_{k, i} d_k,
  \end{equation}
  the following inequality is satisfied:
  \begin{equation}
    \label{Backtrack:Test}
    f(x_k) - f(x_{k + 1, i})
    \geq
    \eta_1 \frac{\DualPairing{\Gradient f(x_k)}{d_k}^2}{L_{k, i} \Norm{d_k}^2}
    \quad (\equiv \eta_1 h_{k, i} \DualPairing{\Gradient f(x_k)}{d_k}),
  \end{equation}
  where $\eta_1 \in \OpenClosedInterval{0}{\frac{1}{2}}$ is a certain (absolute)
  constant.
  Once $i_k$ has been found, we set
  \begin{equation}
    \label{Backtrack:FinalQuants}
    L_k
    \DefinedEqual
    L_{k, i_k},
    \qquad
    h_k
    \DefinedEqual
    h_{k, i_k},
    \qquad
    x_{k + 1}
    \DefinedEqual
    x_{k + 1, i_k}.
  \end{equation}
  For this strategy, we can guarantee that it is well-defined (i.e., $i_k$ is
  finite) and that
  \begin{equation}
    \label{Backtrack:MaxL}
    L_k \leq (1 - \eta_1)^{-1} L,
  \end{equation}
  which means, in view of \cref{Backtrack:Test,Backtrack:FinalQuants}, that
  \cref{FuncProg} is satisfied with
  \[
    \eta \DefinedEqual \eta_1 (1 - \eta_1).
  \]
  Indeed, by construction, $L_{k, i}$ monotonically increases in~$i$.
  Therefore, at some point, it will become greater or equal than~$L$.
  Suppose that $L_{k, i} \geq L$ for some $i \geq 0$.
  Then, from \cref{LipGrad-FuncUbd,Backtrack:Quants}, it follows that
  \[
    \begin{aligned}
      f(x_k) - f(x_{k + 1, i})
      &\geq
      \DualPairing{\Gradient f(x_k)}{x_k - x_{k + 1, i}}
      -
      \frac{L_{k, i}}{2} \Norm{x_k - x_{k + 1, i}}^2
      \\
      &=
      h_{k, i} \DualPairing{\Gradient f(x_k)}{d_k}
      -
      \frac{L_{k, i}}{2} h_{k, i}^2 \Norm{d_k}^2
      \\
      &=
      \frac{
        \DualPairing{\Gradient f(x_k)}{d_k}^2
      }{
        2 L_{k, i} \Norm{d_{k, i}}^2
      }.
    \end{aligned}
  \]
  Thus, inequality~\eqref{Backtrack:Test} will surely be satisfied once
  $L_{k, i}$ becomes greater or equal than~$L$.
  To prove \cref{Backtrack:MaxL} we need to consider two cases.
  If $i_k = 0$, then $L_k = L_{k, 0} \leq L$ according to
  \cref{Backtrack:StartL-Ubd}, so \cref{Backtrack:MaxL} is satisfied.
  If $i_k > 0$, then inequality~\eqref{Backtrack:Test} was not satisfied for
  $i = i_k - 1$, which means that
  \begin{equation}
    \label{Backtrack:Rej}
    f(x_k) - f(x_{k + 1, i_k - 1})
    <
    \eta_1 h_{k, i_k - 1} \DualPairing{\Gradient f(x_k)}{d_k}.
  \end{equation}
  On the other hand, from \cref{LipGrad-FuncUbd,Backtrack:Quants}, it follows
  that
  \[
    \begin{aligned}
      f(x_k) - f(x_{k + 1, i_k - 1})
      &\geq
      \DualPairing{\Gradient f(x_k)}{x_k - x_{k + 1, i_k - 1}}
      -
      \frac{L}{2} \Norm{x_k - x_{k + 1, i_k - 1}}^2
      \\
      &=
      h_{k, i_k - 1} \DualPairing{\Gradient f(x_k)}{d_k}
      -
      \frac{L}{2} h_{k, i_k - 1}^2 \Norm{d_k}^2.
    \end{aligned}
  \]
  Combining this inequality with \cref{Backtrack:Rej}, we obtain
  \[
    h_{k, i_k - 1}
    \geq
    2 (1 - \eta_1)
    \frac{\DualPairing{\Gradient f(x_k)}{d_k}}{L \Norm{d_k}^2}.
  \]
  In view of \cref{Backtrack-Unit:Quants}, it means that
  \[
    L_{k, i_k - 1}
    \geq
    \frac{1}{2 (1 - \eta_1)} L.
  \]
  This proves \cref{Backtrack:MaxL} since $L_k = L_{k, i_k} = 2 L_{k, i_k - 1}$.
\end{example}

\begin{example}[Armijo--Goldstein line search]
  \label{St:Armijo-Goldstein}
  Let $h_k$ be chosen in such a way so that
  \begin{equation}
    \label{Armijo-Goldstein}
    \eta_1 \DualPairing{\Gradient f(x_k)}{x_k - x_{k + 1}}
    \leq f(x_k) - f(x_{k + 1}) \leq
    \eta_2 \DualPairing{\Gradient f(x_k)}{x_k - x_{k + 1}},
  \end{equation}
  where $0 < \eta_1 \leq \eta_2 < 1$ are certain (absolute) constants.
  In this case, \cref{FuncProg} is satisfied with
  \begin{equation}
    \label{Armijo-Goldstein-Prog-Const}
    \eta \DefinedEqual 2 \eta_1 (1 - \eta_2).
  \end{equation}
  Indeed, according to \cref{LipGrad-FuncUbd}, we have
  \[
    f(x_k) - f(x_{k + 1})
    \geq
    \DualPairing{\Gradient f(x_k)}{x_k - x_{k + 1}}
    -
    \frac{L}{2} \Norm{x_k - x_{k + 1}}^2.
  \]
  Combining this with the second inequality in \cref{Armijo-Goldstein}, we
  obtain
  \[
    (1 - \eta_2) \DualPairing{\Gradient f(x_k)}{x_k - x_{k + 1}}
    \leq
    \frac{L}{2} \Norm{x_k - x_{k + 1}}^2.
  \]
  Combining this further with the first inequality in \cref{Armijo-Goldstein},
  we get
  \[
    f(x_k) - f(x_{k + 1})
    \geq
    \eta_1 \DualPairing{\Gradient f(x_k)}{x_k - x_{k + 1}}
    \geq
    2 \eta_1 (1 - \eta_2)
    \frac{
      \DualPairing{\Gradient f(x_k)}{x_k - x_{k + 1}}^2
    }{
      L \Norm{x_k - x_{k + 1}}^2
    }.
  \]
  By the definition of $x_{k + 1}$ in \cref{Alg:BFGS}, this is exactly
  \cref{FuncProg} with $\eta$ given by \cref{Armijo-Goldstein-Prog-Const}.
\end{example}

\begin{example}[Wolfe--Powell line search]
  Let $h_k$ be chosen in such a way so that
  \begin{align}
    \label{WolfeConditions:FuncProg}
    & f(x_k) - f(x_{k + 1})
    \geq
    \eta_1 \DualPairing{\Gradient f(x_k)}{x_k - x_{k + 1}},
    \\
    \label{WolfeConditions:Der}
    & \DualPairing{\Gradient f(x_{k + 1})}{x_k - x_{k + 1}}
    \leq
    \eta_2' \DualPairing{\Gradient f(x_k)}{x_k - x_{k + 1}},
  \end{align}
  where $0 < \eta_1 \leq \eta_2' < 1$ are certain (absolute) constants.
  (The condition $\eta_1 \leq \eta_2'$ guarantees that there exists a step size
  $h_k$ satisfying \cref{WolfeConditions:FuncProg,WolfeConditions:Der}.)
  In this case, \cref{FuncProg} is satisfied with
  \[
    \eta \DefinedEqual \eta_1 (1 - \eta_2').
  \]
  Indeed, according to \cref{WolfeConditions:Der,LipGrad}, we have
  \[
    (1 - \eta_2') \DualPairing{\Gradient f(x_k)}{x_k - x_{k + 1}}
    \leq
    \DualPairing{\Gradient f(x_k) - \Gradient f(x_{k + 1})}{x_k - x_{k + 1}}
    \leq
    L \Norm{x_k - x_{k + 1}}^2.
  \]
  Now we can proceed exactly as in \cref{St:Armijo-Goldstein}.
\end{example}

\begin{remark}
  Note that any step size satisfying the Armijo--Goldstein
  conditions~\eqref{Armijo-Goldstein} also satisfies the Wolfe
  conditions~\eqref{WolfeConditions:FuncProg} and~\eqref{WolfeConditions:Der}
  with $\eta_2' = \eta_2$.
  Indeed, by the convexity of~$f$ and the second inequality in
  \cref{Armijo-Goldstein}, we have
  \[
    \DualPairing{\Gradient f(x_{k + 1})}{x_k - x_{k + 1}}
    \leq
    f(x_k) - f(x_{k + 1})
    \leq
    \eta_2 \DualPairing{\Gradient f(x_k)}{x_k - x_{k + 1}},
  \]
  which is exactly condition~\eqref{WolfeConditions:Der} with
  $\eta_2' = \eta_2$.
\end{remark}

Let us now establish linear convergence of \cref{Alg:BFGS}.
In what follows, we denote
\begin{equation}
  \label{Def-Psi0}
  \psi_0 \DefinedEqual \psi(L B, H_0^{-1}) \ (\geq 0).
\end{equation}
(Recall that $\psi$ is the Bregman distance, defined in \cref{BregmanDist}.)

\begin{lemma}
  \label{St:Bound-Sum-Phi-LB}
  In \cref{Alg:BFGS}, for all $k \geq 1$, we have
  \[
    \sum_{i = 0}^{k - 1}
    \phi_{L B}\bigl( d_i, \Gradient f(x_i) \bigr)
    \leq
    \psi_0
    +
    k \ln \kappa.
  \]
\end{lemma}

\begin{proof}
  Denote $G_k \DefinedEqual H_k^{-1}$ for each $k \geq 0$.
  Then, $G_{k + 1} = \BFGS(G_k, \delta_k, \gamma_k)$ for all $k \geq 0$ and
  $\psi_0 = \psi(L B, G_0)$.

  Let $k \geq 1$ be an arbitrary integer.
  Applying \cref{St:BregmanDist-Upd} with $C \DefinedEqual L B$, we obtain,
  for all $i \geq 0$,
  \[
    \psi(L B, G_i) - \psi(L B, G_{i + 1})
    =
    \phi_{L B}(\delta_i, G \delta_i) - \phi_{L B}(\delta_i, \gamma_i).
  \]
  Summing up these identities for all $i = 0, \dots, k - 1$ and using the fact
  that $\psi(\cdot, \cdot)$ is nonnegative, we get
  \begin{equation}
    \label{Proof:Bound-Sum-Phi-LB:Sum-Prel}
    \sum_{i = 0}^{k - 1}
    \phi_{L B}(\delta_i, G \delta_i)
    \leq
    \psi_0
    +
    \sum_{i = 0}^{k - 1} \phi_{L B}(\delta_i, \gamma_i).
  \end{equation}
  Let us now estimate the sum on the right from above.

  Let $0 \leq i \leq k - 1$ be an arbitrary index.
  Note from \cref{Def-Alpha-Beta} that
  \[
    \begin{aligned}
      &
      \alpha_{L B}(\delta_i, \gamma_i)
      =
      \frac{
        \RelativeDualNorm{\gamma_i}{L B}^2
      }{
        \DualPairing{\gamma_i}{\delta_i}
      }
      =
      \frac{\DualNorm{\gamma_i}^2}{L \DualPairing{\gamma_i}{\delta_i}},
      \\
      &
      \beta_{L B}(\delta_i, \gamma_i)
      =
      \frac{\DualPairing{\gamma_i}{\delta_i}}{\RelativeNorm{\delta_i}{L B}^2}
      =
      \frac{\DualPairing{\gamma_i}{\delta_i}}{L \Norm{\delta_i}^2}.
    \end{aligned}
  \]
  Hence, using the definitions of $\delta_i$ and $\gamma_i$
  (\cref{Alg:BFGS:Def-delta-gamma} of \cref{Alg:BFGS}) and applying
  \cref{LipGrad-ViaGrad,StrongConv,ConditionNumber}, we obtain
  \[
    \alpha_{L B}(\delta_i, \gamma_i)
    \leq
    1,
    \qquad
    \beta_{L B}(\delta_i, \gamma_i)
    \geq
    \kappa^{-1}.
  \]
  Therefore, according to \cref{Def-phi,Def-phi-Short},
  \[
    \phi_{L B}(\delta_i, \gamma_i)
    =
    \alpha_{L B}(\delta_i, \gamma_i) - 1
    -
    \ln \beta_{L B}(\delta_i, \gamma_i)
    \leq
    \ln \kappa.
  \]
  Thus,
  \[
    \sum_{i = 0}^{k - 1} \phi_{L B}(\delta_i, \gamma_i)
    \leq
    k \ln \kappa.
  \]

  Substituting this inequality into \cref{Proof:Bound-Sum-Phi-LB:Sum-Prel}, we
  obtain
  \[
    \sum_{i = 0}^{k - 1} \phi_{L B}(\delta_i, G \delta_i)
    \leq
    \psi_0
    +
    k \ln \kappa.
  \]

  It remains to note that, for all $0 \leq i \leq k - 1$, we have
  \[
    \phi_{L B}(\delta_i, G \delta_i)
    =
    \phi_{L B}(d_i, G d_i)
    =
    \phi_{L B}\bigl( d_i, \Gradient f(x_i) \bigr)
  \]
  since $\delta_i = -h_i d_i$, $d_i = G_i^{-1} \Gradient f(x_i)$
  (see \cref{Alg:BFGS}) and in view of \cref{Alpha-Beta-Phi-Homogeneity}.
\end{proof}

\begin{lemma}
  \label{St:AuxiliaryUbd-Product}
  Let $u_1, \dots, u_k \in \RealField$ and $U \in \RealField$ be such that
  \begin{equation}
    \label{AuxiliaryUbd-Product:Sum-Arg-Ubd}
    \sum_{i = 1}^k u_i \leq U.
  \end{equation}
  Further, let $a \geq 0$ be such that $a e^{-u_i} \leq 1$ for all
  $1 \leq i \leq k$.
  Then, $a e^{-U / k} \leq 1$, and
  \begin{equation}
    \label{AuxiliaryUbd-Product:Result}
    \prod_{i = 1}^k (1 - a e^{-u_i})
    \leq
    (1 - a e^{-U / k})^k.
  \end{equation}
\end{lemma}

\begin{proof}
  Let $\Map{\xi}{\RealField}{\RealField}$ be the function
  \[
    \xi(u) \DefinedEqual a e^{-u}.
  \]
  Note that $\xi$ is a decreasing convex function.
  By our assumptions, $\xi(u_i) \leq 1$ for all $1 \leq i \leq k$.
  From \cref{AuxiliaryUbd-Product:Sum-Arg-Ubd}, we obtain, using the
  monotonicity and convexity of~$\xi$, that
  \begin{equation}
    \label{Proof:AuxiliaryUbd-Product:AuxUbd}
    \xi(U / k)
    \leq
    \xi\Bigl( \frac{1}{k} \sum_{i = 1}^k u_i \Bigr)
    \leq
    \frac{1}{k} \sum_{i = 1}^k \xi(u_i)
    \leq
    1.
  \end{equation}
  This proves that $a e^{-U / k} \leq 1$.

  Now let us prove \cref{AuxiliaryUbd-Product:Result}.
  We can assume that $\xi(u_i) < 1$ for all $1 \leq i \leq k$ since otherwise
  \cref{AuxiliaryUbd-Product:Result} is trivial.
  In particular, from \cref{Proof:AuxiliaryUbd-Product:AuxUbd}, it follows that
  $\xi(U / K) < 1$.
  Let $\Map{\rho}{\EffectiveDomain \rho}{\RealField}$ be the function
  \[
    \rho(u)
    \DefinedEqual
    \ln(1 - a e^{-u})
    =
    \ln\bigl( 1 - \xi(u) \bigr),
    \qquad
    \EffectiveDomain \rho
    \DefinedEqual
    \SetBuilder{u \in \RealField}{\xi(u) < 1}.
  \]
  Note that $\rho$ is an increasing concave function, as the composition of the
  decreasing concave function
  $t \in \ClosedOpenInterval{-\infty}{1} \mapsto \ln(1 - t)$
  and the decreasing convex function~$\xi$.
  Therefore, in view of \cref{AuxiliaryUbd-Product:Sum-Arg-Ubd},
  \[
    \frac{1}{k} \sum_{i = 1}^k \rho(u_i)
    \leq
    \rho\Bigl( \frac{1}{k} \sum_{i = 1}^k u_i \Bigr)
    \leq
    \rho(U / k).
  \]
  Substituting the definition of $\rho$, this inequality reads
  \[
    \frac{1}{k} \sum_{i = 1}^k \ln(1 - a e^{-u_i})
    \leq
    \ln(1 - a e^{-U / k}).
  \]
  Taking the exponents from both sides and rearranging, we obtain
  \cref{AuxiliaryUbd-Product:Result}.
\end{proof}

\begin{theorem}
  \label{St:LinearConvergence}
  In \cref{Alg:BFGS}, for all $k \geq 1$, we have
  \begin{equation}
    \label{LinearConvergence}
    f(x_k) - f_*
    \leq
    (1 - 2 \eta \kappa^{-2} \zeta_k)^k [f(x_0) - f_*],
  \end{equation}
  where $\zeta_k \DefinedEqual \exp(-\psi_0 / k)$.
\end{theorem}

\begin{proof}
  \ProofPart

  Let $k \geq 0$ be arbitrary.
  Denote
  \begin{equation}
    \label{Proof:LinearConvergence:Def-theta}
    \theta_k
    \DefinedEqual
    \frac{
      \DualPairing{\Gradient f(x_k)}{d_k}^2
    }{
      \DualNorm{\Gradient f(x_k)}^2 \Norm{d_k}^2
    }
    \quad (\leq 1).
  \end{equation}
  Then, in view of \cref{As-StepSize-Global,StrongConv-FuncUbd}, we have
  \[
    f(x_k) - f(x_{k + 1})
    \geq
    \eta \frac{\DualPairing{\Gradient f(x_k)}{d_k}^2}{L \Norm{d_k}^2}
    =
    \eta \theta_k \frac{1}{L} \DualNorm{\Gradient f(x_k)}^2
    \geq
    2 \eta \theta_k \frac{\mu}{L} [f(x_k) - f_*].
  \]
  In other words, according to \cref{ConditionNumber},
  \[
    f(x_{k + 1}) - f^*
    \leq
    (1 - 2 \eta \kappa^{-1} \theta_k) [f(x_k) - f^*].
  \]

  In particular, for all $k \geq 1$, it follows that
  \begin{equation}
    \label{Proof:LinearConvergence:Prel}
    f(x_k) - f^*
    \leq
    [f(x_0) - f^*] \prod_{i = 0}^{k - 1} (1 - 2 \eta \kappa^{-1} \theta_i),
  \end{equation}
  where
  \begin{equation}
    \label{Proof:LinearConvergence:Coef}
    2 \eta \kappa^{-1} \theta_i \leq 1
  \end{equation}
  for all $0 \leq i \leq k - 1$.
  Let us estimate the latter product from above.

  \ProofPart

  Let $k \geq 1$ be an arbitrary index.
  According to \cref{St:Bound-Sum-Phi-LB}, we have
  \begin{equation}
    \label{Proof:LinearConvergence:Sum-Phi-Prel}
    \sum_{i = 0}^{k - 1}
    \phi_{L B}\bigl( d_i, \Gradient f(x_i) \bigr)
    \leq
    \psi_0 + k \ln \kappa.
  \end{equation}
  Let $0 \leq i \leq k - 1$ be an arbitrary index.
  In view of \cref{Def-phi}, for any $\alpha \geq \beta > 0$, we have
  \[
    \phi(\alpha, \beta)
    =
    \alpha - 1 - \ln \beta
    =
    \alpha - 1 - \ln \alpha
    +
    \ln \frac{\alpha}{\beta}
    \geq
    \ln \frac{\alpha}{\beta}.
  \]
  Thus, according to
  \cref{Def-phi-Short,Def-Alpha-Beta,Proof:LinearConvergence:Def-theta},
  \[
    \phi_{L B}\bigl( d_i, \Gradient f(x_i) \bigr)
    \geq
    \ln \frac{
      \alpha_{L B}\bigl(d_i, \Gradient f(x_i) \bigr)
    }{
      \beta_{L B}\bigl( d_i, \Gradient f(x_i) \bigr)
    }
    =
    \ln \frac{
      \DualNorm{\Gradient f(x_i)}^2 \Norm{d_i}^2
    }{
      \DualPairing{\Gradient f(x_i)}{d_i}^2
    }
    =
    -\ln \theta_i
    \EqualDefines
    u_i.
  \]
  Combining this with \cref{Proof:LinearConvergence:Sum-Phi-Prel}, we obtain
  \[
    \sum_{i = 0}^{k - 1} u_i
    \leq
    \psi_0 + k \ln \kappa.
  \]
  Applying \cref{St:AuxiliaryUbd-Product}
  (using \cref{Proof:LinearConvergence:Coef}), we obtain
  \[
    \prod_{i = 0}^{k - 1} (1 - 2 \eta \kappa^{-1} e^{-u_i})
    \leq
    (1 - 2 \eta \kappa^{-1} e^{-\psi_0 / k + \ln \kappa})^k
    =
    (1 - 2 \eta \kappa^{-2} \zeta_k)^k.
  \]
  Substituting this into \cref{Proof:LinearConvergence:Prel}, we obtain
  \cref{LinearConvergence}.
\end{proof}

Let us discuss the result of \cref{St:LinearConvergence}.
We see that the rate of linear convergence of \cref{Alg:BFGS} is affected by
the choice of the initial inverse Hessian approximation~$H_0$.
The best choice corresponds to
\begin{equation}
  \label{InitialHessApprox-LB}
  H_0 = (L B)^{-1},
\end{equation}
in which case we get $\zeta_k = 1$ for all $k \geq 1$, and
\[
  f(x_k) - f^* \leq (1 - 2 \eta \kappa^{-2})^k [f(x_0) - f^*].
\]
This is the global linear convergence with constant $2 \eta \kappa^{-2}$.
Nevertheless, even if we use any other initial Hessian approximation instead of
\cref{InitialHessApprox-LB}, we can still get almost the same convergence
result.
The only difference is that this linear convergence might begin after some
initial number of iterations.

Indeed, let us introduce the following
\emph{starting moment of global linear convergence}:
\begin{equation}
  \label{StartingMoment-Global}
  K_0^G \DefinedEqual \Ceil{2 \psi_0} \ (\geq 0).
\end{equation}
(Here ``G'' stands for ``global'').
Then, we can prove the following result.

\begin{corollary}
  \label{St:LinearConvergence-Started}
  For all $k \geq K_0^G$, we have
  \begin{equation}
    \label{LinearConvergence-Started}
    f(x_k) - f^*
    \leq
    (1 - \eta \kappa^{-2})^k [f(x_0) - f^*].
  \end{equation}
\end{corollary}

\begin{proof}
  Indeed, for any $k \geq K_0^G$, according to \cref{StartingMoment-Global}, we
  have $k \geq 2 \psi_0$, and hence
  \[
    \zeta_k
    \equiv
    \exp(-\psi_0 / k)
    \geq
    \exp(-1 / 2)
    \geq
    1 - \frac{1}{2}
    =
    \frac{1}{2}.
  \]
  (We assume that $k \neq 0$ since otherwise \cref{LinearConvergence-Started} is
  trivial.)
  Applying \cref{St:LinearConvergence}, we get \cref{LinearConvergence-Started}.
\end{proof}

Note that, for any ``reasonable'' initial Hessian approximation, namely, for
any~$H_0$ satisfying
\begin{equation}
  \label{ResonableInitialHessApprox}
  (L B)^{-1} \preceq H_0 \preceq (\mu B)^{-1},
\end{equation}
we have a ``reasonable'' starting moment of global linear convergence:
\[
  K_0^G \leq [2 n \ln \kappa] \quad (\ll \kappa^2).
\]
Indeed, denote $G_0 \DefinedEqual H_0^{-1}$.
Then,
\[
  \mu B \preceq G_0 \preceq L B,
\]
which means that each eigenvalue~$\lambda_i$ of~$G_0$ w.r.t.~$L B$ is bounded
as follows:
\[
  \kappa^{-1} \leq \lambda_i \leq 1.
\]
Hence, according to \cref{Def-Psi0,BregmanDist-ViaEigs,Def-omega} and the fact
that $\omega$ monotonically decreases on $\OpenClosedInterval{-1}{0}$, we have
\[
  \psi_0
  =
  \sum_{i = 1}^n \omega(\lambda_i - 1)
  \leq
  n \omega(\kappa^{-1} - 1)
  =
  n [\kappa^{-1} - 1 - \ln \kappa^{-1}]
  \leq
  n \ln \kappa.
\]
Note that it is very easy to generate $H_0$ satisfying
\cref{ResonableInitialHessApprox} even without the knowledge of constants~$\mu$
and~$L$.
One simple strategy is to choose any two points $x_0 \neq x_0'$, compute
\[
  \delta_0' \DefinedEqual x_0' - x_0,
  \qquad
  \gamma_0' \DefinedEqual \Gradient f(x_0') - \Gradient f(x_0),
\]
and set
\[
  H_0
  \DefinedEqual
  \frac{
    \DualPairing{\gamma_0'}{\delta_0'}
  }{
    \DualNorm{\gamma_0'}^2
  } B^{-1}
  \qquad \text{or} \qquad
  H_0
  \DefinedEqual
  \frac{
    \Norm{\delta_0'}^2
  }{
    \DualPairing{\gamma_0'}{\delta_0'}
  } B^{-1}.
\]
These are the standard strategies commonly used for computing the initial
Hessian approximation in the BFGS Method (see, e.g., Section~6.1
in~\cite{Nocedal.Wright-NumericalOptimization-06}).

\subsection{Backtracking Line Search Starting From Unit Step Size}

Let us present the analysis of the backtracking line search strategy, similar
to the one from \cref{St:Backtrack}, but which tries the unit step size
$h_{k, 0} \equiv 1$ first.
Specifically, we find the smallest index $i_k \geq 0$ over
$i = 0, 1, 2, \ldots$ such that, for
\begin{equation}
  \label{Backtrack-Unit:Quants}
  h_{k, i} \DefinedEqual 2^{-i},
  \qquad
  x_{k + 1, i} \DefinedEqual x_k - h_{k, i} d_k,
\end{equation}
the following inequality is satisfied:
\begin{equation}
  \label{Backtrack-Unit:Test}
  f(x_k) - f(x_{k + 1, i})
  \geq
  \eta_1 h_{k, i} \DualPairing{\Gradient f(x_k)}{d_k},
\end{equation}
where $\eta_1 \in \OpenClosedInterval{0}{\frac{1}{2}}$ is a certain (absolute)
constant.
Once $i_k$ has been found, we set
\begin{equation}
  \label{Backtrack-Unit:FinalQuants}
  h_k \DefinedEqual h_{k, i_k},
  \qquad
  x_{k + 1} \DefinedEqual x_{k + 1, i_k}.
\end{equation}

Note that this is exactly the strategy from \cref{St:Backtrack} with the choice
of $L_{k, 0}$ which guarantees that $h_{k, 0} = 1$, i.e.,
\[
  L_{k, 0}
  \DefinedEqual
  \frac{\DualPairing{\Gradient f(x_k)}{d_k}}{\Norm{d_k}^2}.
\]
However, we cannot analyze it using \cref{As-StepSize-Global} since we cannot
guarantee that condition~\eqref{Backtrack:StartL-Ubd} is satisfied.

Nevertheless, we can slightly modify our previous analysis and obtain some the
same result as in \cref{St:LinearConvergence}, up to an absolute constant.

First, let us prove a useful auxiliary inequality.

\begin{lemma}
  \label{St:Omega-Lbd-Via-Log}
  For any $\alpha > 0$, we have
  \[
    \omega(\alpha - 1) \geq \ln \frac{\alpha}{2}.
  \]
\end{lemma}

\begin{proof}
  By \cref{Def-omega}, we need to prove that
  \[
    \alpha - 1 - \ln \alpha \geq \ln \frac{\alpha}{2},
  \]
  or, equivalently, that
  \[
    2 \ln \alpha \leq \alpha - 1 + \ln 2.
  \]
  But this easily follows from the concavity of the logarithm: indeed, we have
  \[
    \ln \alpha \leq \ln 2 + \frac{1}{2} (\alpha - 2).
  \]
  Hence,
  \[
    2 \ln \alpha \leq 2 \ln 2 + \alpha - 2,
  \]
  and we need to check if
  \[
    2 \ln 2 + \alpha - 2 \leq \alpha - 1 + \ln 2.
  \]
  But this is indeed true since $\ln 2 \leq 1$.
\end{proof}

\begin{theorem}
  \label{St:LinearConvergence-Backtrack-Unit}
  Consider \cref{Alg:BFGS} with the backtracking line
  search~\eqref{Backtrack-Unit:Quants}, \eqref{Backtrack-Unit:Test}
  and~\eqref{Backtrack-Unit:FinalQuants}.
  Then, for all $k \geq 1$, we have
  \[
    f(x_k) - f_*
    \leq
    (1 - \eta_1 \kappa^{-2} \zeta_k)^k [f(x_0) - f^*].
  \]
\end{theorem}

\begin{proof}
  \ProofPart

  Let $k \geq 0$ be an arbitrary index.
  If $i_k = 0$, then $h_k = 1$ and so we have, according to
  \cref{Backtrack-Unit:Test},
  \[
    f(x_k) - f(x_{k + 1})
    \geq
    \eta_1 \DualPairing{\Gradient f(x_k)}{d_k}
    =
    \eta_1 \frac{1}{L} \alpha_k^{-1} \DualNorm{\Gradient f(x_k)}^2,
  \]
  where
  \begin{equation}
    \label{Proof:LinearConvergence-Backtrack-Unit:Def-alpha}
    \alpha_k
    \DefinedEqual
    \alpha_{L B}\bigl( d_k, \Gradient f(x_k) \bigr)
    =
    \frac{
      \DualNorm{\Gradient f(x_k)}^2
    }{
      L \DualPairing{\Gradient f(x_k)}{d_k}
    }
  \end{equation}
  (see \cref{Def-Alpha-Beta}).
  If $i_k > 0$, then $h_k \leq \frac{1}{2}$ and we have, according to
  \cref{Backtrack-Unit:Test},
  \begin{equation}
    \label{Proof:LinearConvergence-Backtrack-Unit:StepSizeAcc}
    f(x_k) - f(x_{k + 1})
    \geq
    \eta_1 h_k \DualPairing{\Gradient f(x_k)}{d_k}.
  \end{equation}
  At the same time, since the step size $h_{k, i_k - 1}$ was not accepted, we
  have
  \begin{equation}
    \label{Proof:LinearConvergence-Backtrack-Unit:StepSizeRej}
    f(x_k) - f(x_{k, i_k - 1})
    <
    \eta_1 h_{k, i_k - 1} \DualPairing{\Gradient f(x_k)}{d_k}.
  \end{equation}
  On the other hand, in view of \cref{LipGrad-FuncUbd,Backtrack-Unit:Quants},
  \[
    \begin{aligned}
      f(x_k) - f(x_{k, i_k - 1})
      &\geq
      \DualPairing{\Gradient f(x_k)}{h_{k, i_k - 1} d_k}
      -
      \frac{L}{2} \Norm{h_{k, i_k - 1} d_k}^2
      \\
      &=
      h_{k, i_k - 1} \DualPairing{\Gradient f(x_k)}{d_k}
      -
      \frac{L}{2} h_{k, i_k - 1}^2 \Norm{d_k}^2.
    \end{aligned}
  \]
  Combining this with \cref{Proof:LinearConvergence-Backtrack-Unit:StepSizeRej},
  we obtain
  \[
    h_{k, i_k - 1}
    \geq
    2 (1 - \eta_1) \frac{\DualPairing{\Gradient f(x_k)}{d_k}}{L \Norm{d_k}^2}.
  \]
  Therefore, in view of \cref{Backtrack-Unit:FinalQuants,Backtrack-Unit:Quants},
  \[
    h_k
    \equiv
    h_{k, i_k}
    =
    \frac{1}{2} h_{k, i_k - 1}
    \geq
    (1 - \eta_1) \frac{\DualPairing{\Gradient f(x_k)}{d_k}}{L \Norm{d_k}^2}.
  \]
  Substituting this inequality into
  \cref{Proof:LinearConvergence-Backtrack-Unit:StepSizeAcc}, we conclude that,
  in the case when $i_k > 0$,
  \[
    f(x_k) - f(x_{k + 1})
    \geq
    \eta_1 (1 - \eta_1)
    \frac{\DualPairing{\Gradient f(x_k)}{d_k}^2}{L \Norm{d_k}^2}
    =
    \eta_1 (1 - \eta_1) \frac{1}{L} \theta_k \DualNorm{\Gradient f(x_k)}^2,
  \]
  where
  \begin{equation}
    \label{Proof:LinearConvergence-Backtrack-Unit:Def-theta}
    \theta_k
    \DefinedEqual
    \frac{
      \DualPairing{\Gradient f(x_k)}{d_k}^2
    }{
      \DualNorm{\Gradient f(x_k)}^2 \Norm{d_k}^2
    }
    \quad (\in \OpenClosedInterval{0}{1}).
  \end{equation}

  \ProofPart

  Thus, we obtain:
  \[
    f(x_k) - f(x_{k + 1})
    \geq
    \begin{cases}
      \eta_1 \alpha_k^{-1} \frac{1}{L} \DualNorm{\Gradient f(x_k)}^2,
      & \text{if $h_k = 1$},
      \\
      \eta_1 (1 - \eta_1) \theta_k \frac{1}{L} \DualNorm{\Gradient f(x_k)}^2,
      & \text{otherwise}.
    \end{cases}
  \]
  Denote
  \[
    I_1
    \DefinedEqual
    \SetBuilder{0 \leq i \leq k - 1}{h_i = 1},
    \qquad
    \bar{I}_1
    \DefinedEqual
    \SetBuilder{0 \leq i \leq k - 1}{i \notin I_1}.
  \]
  Then,
  \[
    f(x_i) - f(x_{i + 1})
    \geq
    \begin{cases}
      \eta_1 \alpha_i^{-1} \frac{1}{L} \DualNorm{\Gradient f(x_i)}^2,
      & \text{if $i \in I_1$},
      \\
      \eta_1 (1 - \eta_1) \theta_i \frac{1}{L} \DualNorm{\Gradient f(x_i)}^2,
      & \text{if $i \in \bar{I}_1$}.
    \end{cases}
  \]
  Applying \cref{StrongConv-FuncUbd}, we obtain
  $\DualNorm{\Gradient f(x_i)}^2 \geq 2 \mu [f(x_i) - f^*]$ and thus,
  in view of \cref{ConditionNumber},
  \[
    f(x_i) - f(x_{i + 1})
    \geq
    \begin{cases}
      2 \eta_1 \kappa^{-1} \alpha_i^{-1} [f(x_i) - f^*],
      & \text{if $i \in I_1$},
      \\
      2 \eta_1 (1 - \eta_1) \kappa^{-1} \theta_i [f(x_i) - f^*],
      & \text{if $i \in \bar{I}_1$}.
    \end{cases}
  \]
  For any $0 \leq i \leq k - 1$, denote
  \begin{equation}
    \label{Proof:LinearConvergence-Backtrack-Unit:Def-xi}
    \hat{\theta}_i
    \DefinedEqual
    \begin{cases}
      2 / \alpha_i, & \text{if $i \in I_1$},
      \\
      \theta_i, & \text{if $i \in \bar{I}_1$}.
    \end{cases}
  \end{equation}
  Since $2 (1 - \eta_1) \geq 1$ (since $\eta_1 \leq \frac{1}{2}$), we thus have,
  for all $0 \leq i \leq k - 1$,
  \[
    f(x_i) - f(x_{i + 1})
    \geq
    \eta_1 \kappa^{-1} \hat{\theta}_i [f(x_i) - f^*],
  \]
  or, equivalently,
  \[
    f(x_{i + 1}) - f^*
    \leq
    (1 - \eta_1 \kappa^{-1} \hat{\theta}_i) [f(x_i) - f^*].
  \]
  It follows that, for all $k \geq 1$, we have
  \begin{equation}
    \label{Proof:LinearConvergence-Backtrack-Unit:Prel}
    f(x_k) - f^*
    \leq
    [f(x_0) - f^*]
    \prod_{i = 0}^{k - 1} (1 - \eta_1 \kappa^{-1} \hat{\theta}_i),
  \end{equation}
  where
  \begin{equation}
    \label{Proof:LinearConvergence-Backtrack-Unit:Coef}
    \eta_1 \kappa^{-1} \hat{\theta}_i \leq 1
  \end{equation}
  for all $0 \leq i \leq k - 1$.
  It remains to estimate the product from above.

  \ProofPart

  According to \cref{St:Bound-Sum-Phi-LB}, we have
  \[
    \sum_{i = 0}^{k - 1} \phi_{L B}\bigl( d_i, \Gradient f(x_i) \bigr)
    \leq
    \psi_0 + k \ln \kappa.
  \]
  Note that, in view of \cref{Def-phi,Def-omega},
  for any $\alpha \geq \beta > 0$, we have
  \[
    \phi(\alpha, \beta)
    =
    \alpha - 1 - \ln \beta
    =
    \alpha - 1 - \ln \alpha - \ln \frac{\beta}{\alpha}
    =
    \omega(\alpha - 1) - \ln \frac{\beta}{\alpha}.
  \]
  From \cref{Def-Alpha-Beta,Proof:LinearConvergence-Backtrack-Unit:Def-theta},
  note that
  \[
    \frac{
      \beta_{L B}\bigl( d_i, \Gradient f(x_i) \bigr)
    }{
      \alpha_{L B}\bigl( d_i, \Gradient f(x_i) \bigr)
    }
    =
    \frac{
      \DualPairing{\Gradient f(x_i)}{d_i}^2
    }{
      \DualNorm{\Gradient f(x_i)}^2 \Norm{d_i}^2
    }
    =
    \theta_i.
  \]
  Thus, in view of
  \cref{Proof:LinearConvergence-Backtrack-Unit:Def-alpha,Def-phi-Short},
  \[
    \phi_{L B}\bigl( d_i, \Gradient f(x_i) \bigr)
    =
    \omega(\alpha_i - 1) - \ln \theta_i.
  \]
  Using the fact that, for any $0 \leq i \leq k - 1$, we have
  $\omega(\alpha_i - 1) \geq 0$ and $-\ln \theta_i \geq 0$
  (since $\theta_i \leq 1$), we obtain
  \[
    \begin{aligned}
      \hspace{2em}&\hspace{-2em}
      \sum_{i = 0}^{k - 1} \phi_{L B}\bigl( d_i, \Gradient f(x_i) \bigr)
      =
      \sum_{i = 0}^{k - 1} [\omega(\alpha_i - 1) - \ln \theta_i]
      \\
      &=
      \sum_{i \in I_1} [\omega(\alpha_i - 1) - \ln \theta_i]
      +
      \sum_{i \in \bar{I}_1} [\omega(\alpha_i - 1) - \ln \theta_i]
      \\
      &\geq
      \sum_{i \in I_1} \omega(\alpha_i - 1)
      +
      \sum_{i \in \bar{I}_1} -\ln \theta_i.
    \end{aligned}
  \]
  Applying \cref{St:Omega-Lbd-Via-Log} and using
  \cref{Proof:LinearConvergence-Backtrack-Unit:Def-xi}, we obtain,
  for all $i \in I_1$,
  \[
    \omega(\alpha_i - 1)
    \geq
    \ln \frac{\alpha_i}{2}
    =
    -\ln \hat{\theta}_i.
  \]
  At the same time, for all $i \in \bar{I}_1$, in view of
  \cref{Proof:LinearConvergence-Backtrack-Unit:Def-xi},
  \[
    -\ln \theta_i = -\ln \hat{\theta}_i.
  \]
  Thus, we have
  \[
    \sum_{i = 0}^{k - 1} -\ln \hat{\theta}_i
    \leq
    \sum_{i = 0}^{k - 1} \phi_{L B}\bigl( d_i, \Gradient f(x_i) \bigr)
    \leq
    \psi_0 + k \ln \kappa.
  \]
  Applying \cref{St:AuxiliaryUbd-Product}
  (using \cref{Proof:LinearConvergence-Backtrack-Unit:Coef}), we obtain
  \[
    \prod_{i = 0}^{k - 1} (1 - \eta_1 \kappa^{-1} \hat{\theta}_i)
    \leq
    (1 - \eta_1 \kappa^{-2} \zeta_k)^k.
  \]
  It remains to substitute this bound into
  \cref{Proof:LinearConvergence-Backtrack-Unit:Prel}.
\end{proof}

  \section{Local Superlinear Convergence}

Now we additionally assume that the objective function in problem~\eqref{Prob}
has Lipschitz continuous Hessian: there exists a constant $L_2 \geq 0$ such
that, for all $x, y \in \VectorSpace{E}$, it holds
\begin{equation}
  \label{LipHess}
  \Norm{\Hessian f(x) - \Hessian f(y)} \leq L_2 \Norm{x - y}.
\end{equation}

From \cref{LipHess}, it follows that (see Lemma~1.2.4
in~\cite{Nesterov-LecturesConvex-18}), for all $x, y \in \VectorSpace{E}$,
\begin{equation}
  \label{LipHess-FuncBnd}
  \Abs[\Big]{
    f(y)
    -
    f(x)
    -
    \DualPairing{\Gradient f(x)}{y - x}
    -
    \frac{1}{2} \DualPairing{\Hessian f(x) (y - x)}{y - x}
  }
  \leq
  \frac{L_2}{6} \Norm{y - x}^3.
\end{equation}

For functions, satisfying \cref{StrongConv,LipHess}, we have the following
convenient inequalities.
In what follows, we denote
\begin{equation}
  \label{Def-H}
  H
  \DefinedEqual
  \sqrt{2} \frac{L_2}{\mu^{3 / 2}}.
\end{equation}

\begin{lemma}
  \label{St:AuxiliaryBounds-Hess}
  \begin{EnumerateClaims}
    \item
      \label{St:AuxiliaryBounds-Hess:Simple}
      For any $x, y \in \VectorSpace{E}$, we have
      \begin{equation}
        \label{AuxiliaryBounds-Hess:Simple}
        [1 + H \sigma(x, y)]^{-1} \Hessian f(x)
        \preceq
        \Hessian f(y)
        \preceq
        [1 + H \sigma(x, y)] \Hessian f(x),
      \end{equation}
      where
      $
        \sigma(x, y)
        \DefinedEqual
        [f(y) - f(x) - \DualPairing{\Gradient f(x)}{y - x}]^{1 / 2}
      $.

    \item
      \label{St:AuxiliaryBounds-Hess:Int}
      For any $x, y, z \in \VectorSpace{E}$ and
      $J(x, y) \DefinedEqual \int_0^1 \Hessian f\bigl( x + t (y - x) \bigr) dt$,
      we have
      \begin{equation}
        \label{AuxiliaryBounds-Hess:Int}
        [1 + H \sigma(z, x, y)]^{-1} \Hessian f(z)
        \preceq
        J(x, y)
        \preceq
        [1 + H \sigma(z, x, y)] \Hessian f(z),
      \end{equation}
      where
      $
        \sigma(z, x, y)
        \DefinedEqual
        \frac{1}{2} \sigma(z, x) + \frac{1}{2} \sigma(z, y)
      $.

    \item
      \label{St:AuxiliaryBounds-Hess:Int2}
      For any $x, y, z \in \VectorSpace{E}$ and
      $J_2(x, y) \DefinedEqual 2 \int_0^1 (1 - t) \Hessian f(x + t (y - x)) dt$,
      we have
      \begin{equation}
        \label{AuxiliaryBounds-Hess:Int2}
        [1 + H \sigma_2(z, x, y)]^{-1} \Hessian f(z)
        \preceq
        J_2(x, y)
        \preceq
        [1 + H \sigma_2(z, x, y)] \Hessian f(z),
      \end{equation}
      where
      $
        \sigma_2(z, x, y)
        \DefinedEqual
        \frac{2}{3} \sigma(z, x) + \frac{1}{3} \sigma(z, y)
      $.
  \end{EnumerateClaims}
\end{lemma}

\begin{proof}
  Note that \cref{LipHess,StrongConv} are equivalent to the following operator
  inequalities for all $x, y \in \VectorSpace{E}$:
  \begin{gather}
    \label{LipHess-OpIneq}
    \Hessian f(y) - \Hessian f(x)
    \preceq
    L_2 \Norm{y - x} B,
    \\
    \label{StrongConv-OpIneq}
    \Hessian f(x) \succeq \mu B.
  \end{gather}
  Note also that, in view of \cref{StrongConv-FuncLbd}, we have
  \begin{equation}
    \label{StrongConv-FuncLbd-ViaSigma}
    \Norm{y - x} \leq \sqrt{\frac{2}{\mu}} \, \sigma(x, y).
  \end{equation}

  Let us prove \cref{AuxiliaryBounds-Hess:Simple}.
  From \cref{LipHess-OpIneq,StrongConv-OpIneq}, it follows that
  \begin{gather*}
    \Hessian f(y) - \Hessian f(x)
    \preceq
    \frac{L_2}{\mu} \Norm{y - x} \Hessian f(x),
    \\
    \Hessian f(x) - \Hessian f(y)
    \preceq
    \frac{L_2}{\mu} \Norm{y - x} \Hessian f(y).
  \end{gather*}
  Rearranging, we obtain
  \[
    \Bigl( 1 + \frac{L_2}{\mu} \Norm{y - x} \Bigr)^{-1} \Hessian f(x)
    \preceq
    \Hessian f(y)
    \preceq
    \Bigl( 1 + \frac{L_2}{\mu} \Norm{y - x} \Bigr) \Hessian f(x).
  \]
  Applying \cref{StrongConv-FuncLbd-ViaSigma}, we arrive at
  \cref{AuxiliaryBounds-Hess:Simple}.

  Now let us prove \cref{AuxiliaryBounds-Hess:Int}.
  Note that
  \[
    J(x, y) - \Hessian f(z)
    =
    \int_0^1 \bigl[
      \Hessian f\bigl(x + t (y - x) \bigr) - \Hessian f(z)
    \bigr] dt.
  \]
  Applying \cref{LipHess}, we obtain
  \[
    \begin{aligned}
      \Norm{J(x, y) - \Hessian f(z)}
      &\leq
      L_2 \int_0^1 \Norm{(1 - t) (x - z) + t (y - z)} d t
      \\
      &\leq
      L_2 r(x, y, z),
    \end{aligned}
  \]
  where
  \[
    r(x, y, z)
    \DefinedEqual
    \frac{1}{2} (\Norm{x - z} + \Norm{y - z}).
  \]
  Thus,
  \[
    -L_2 r(x, y, z) B
    \preceq
    J(x, y) - \Hessian f(z)
    \preceq
    L_2 r(x, y, z) B
  \]
  Hence, in view of \cref{StrongConv-OpIneq},
  \[
    -\frac{L_2}{\mu} r(x, y, z) \Hessian J(x, y)
    \preceq
    J(x, y) - \Hessian f(z)
    \preceq
    \frac{L_2}{\mu} r(x, y, z) \Hessian f(z).
  \]
  Rearranging, we obtain
  \[
    \Bigl( 1 + \frac{L_2}{\mu} r(x, y, z) \Bigr)^{-1} \Hessian f(z)
    \preceq
    J(x, y)
    \preceq
    \Bigl( 1 + \frac{L_2}{\mu} r(x, y, z) \Bigr) \Hessian f(z).
  \]
  Applying \cref{StrongConv-FuncLbd-ViaSigma}, we obtain
  \cref{AuxiliaryBounds-Hess:Int}.

  Finally, let us prove \cref{AuxiliaryBounds-Hess:Int2}.
  Note that
  \[
    J_2(x, y) - \Hessian f(z)
    =
    2 \int_0^1 (1 - t)
    [\Hessian f\bigl( x + t (y - x) \bigr) - \Hessian f(z)] d t.
  \]
  Hence, in view of \cref{LipHess},
  \[
    \begin{aligned}
      \Norm{J_2(x, y) - \Hessian f(z)}
      &\leq
      2 L_2 \int_0^1 (1 - t) \Norm{(1 - t) (x - z) + t (y - z)} d t
      \\
      &\leq
      2 L_2 \Bigl(
        \Norm{x - z} \int_0^1 (1 - t)^2 d t
        +
        \Norm{y - z} \int_0^1 t (1 - t) d t
      \Bigr)
      \\
      &=
      L_2 r_2(z, x, y),
    \end{aligned}
  \]
  where
  \[
    r_2(z, x, y)
    \DefinedEqual
    \frac{2}{3} \Norm{x - z}
    +
    \frac{1}{3} \Norm{y - z}.
  \]
  The rest of the proof is the same as before.
\end{proof}

\begin{lemma}
  \label{St:LocalMetrics}
  Let $f$ be a $\mu$-strongly convex function with $L_2$-Lipschitz continuous
  Hessian.
  Let $x_*$ be the minimizer of $f$, and let $x \in \VectorSpace{E}$.
  Denote
  \[
    \delta \DefinedEqual f(x) - f_*,
    \qquad
    r_* \DefinedEqual \RelativeNorm{x - x_*}{x_*},
    \qquad
    \lambda \DefinedEqual \RelativeDualNorm{\Gradient f(x)}{x},
    \qquad
    \lambda_* \DefinedEqual \RelativeDualNorm{\Gradient f(x)}{x_*}.
  \]
  Then, we have
  \begin{gather}
    \label{LocalMetrics:Dist-Via-Func}
    (1 + \tfrac{1}{3} H \sqrt{\delta} \,)^{-1} \, 2 \delta
    \leq
    r_*^2
    \leq
    (1 + \tfrac{1}{3} H \sqrt{\delta} \,) \, 2 \delta,
    \\
    \label{LocalMetrics:GradNorm-Via-Func}
    (1 + \tfrac{1}{2} H \sqrt{\delta} \,)^{-2} \, 2 \delta
    \leq
    \lambda^2
    \leq
    (1 + \tfrac{1}{4} H \sqrt{\delta} \,)^2 \, 2 \delta,
    \\
    \label{LocalMetrics:GradNormOpt-Via-Func}
    (1 + \tfrac{1}{6} H \sqrt{\delta} \,)^{-1}
    (1 + \tfrac{1}{2} H \sqrt{\delta} \,)^{-1} \,
    2 \delta
    \leq
    \lambda_*^2
    \leq
    (1 + \tfrac{1}{6} H \sqrt{\delta} \,)
    (1 + \tfrac{1}{2} H \sqrt{\delta} \,) \,
    2 \delta.
  \end{gather}
\end{lemma}

\begin{remark}
  The lower bound in \cref{LocalMetrics:GradNorm-Via-Func} can be improved up to
  \[
    (1 + \tfrac{1}{6} H \sqrt{\delta} \,)^{-1}
    (1 + \tfrac{1}{2} H \sqrt{\delta} \,)^{-1} \,
    2 \delta
    \leq
    \lambda^2.
  \]
  The proof is the same as for the lower bound in
  \cref{LocalMetrics:GradNormOpt-Via-Func}.
\end{remark}

\begin{proof}
  \ProofPart

  First, let us prove \cref{LocalMetrics:Dist-Via-Func}.
  Note that
  \[
    \delta
    =
    \frac{1}{2} \DualPairing{J_2 h}{h},
    \qquad
    h \DefinedEqual x - x_*,
    \qquad
    J_2
    \DefinedEqual
    2 \int_0^1 (1 - t) \Hessian f\bigl( x_* + t h \bigr) d t.
  \]
  By \cref{St:AuxiliaryBounds-Hess:Int2}, we have
  \[
    (1 + \tfrac{1}{3} H \sqrt{\delta} \,)^{-1} \Hessian f(x_*)
    \preceq
    J_2
    \preceq
    (1 + \tfrac{1}{3} H \sqrt{\delta} \,) \Hessian f(x_*).
  \]
  Therefore,
  \[
    (1 + \tfrac{1}{3} H \sqrt{\delta} \,)^{-1} r_*^2
    \leq
    2 \delta
    \leq
    (1 + \tfrac{1}{3} H \sqrt{\delta} \,) r_*^2.
  \]
  Rearranging, we obtain \cref{LocalMetrics:Dist-Via-Func}.

  \ProofPart

  Now let us prove \cref{LocalMetrics:GradNorm-Via-Func}.
  Note that we can treat the function $f$ as self-concordant with parameter
  \[
    M \DefinedEqual \frac{L_2}{2 \mu^{3 / 2}} = \frac{1}{2 \sqrt{2}} H.
  \]
  Then, by Theorem~5.2.1 in~\cite{Nesterov-LecturesConvex-18}, we have the
  following inequalities:
  \[
    \omega(M \lambda)
    \leq
    M^2 \delta
    \leq
    \omega_*(M \lambda),
  \]
  where $\Map{\omega_*}{\OpenClosedInterval{-1}{0}}{\RealField}$ is the function
  $\omega_*(\tau) \DefinedEqual \omega(-\tau) = -\tau - \ln(1 - \tau)$, and the
  upper bound is valid under the additional assumption that $M \lambda < 1$.
  Since both $\omega$ and $\omega_*$ are increasing functions for nonnegative
  arguments, it follows that
  \[
    \omega_*^{-1}(M^2 \delta)
    \leq
    M \lambda
    \leq
    \omega^{-1}(M^2 \delta),
  \]
  where $\omega$ and $\omega_*$ are the corresponding inverse function (for
  $\omega$, we consider only the right branch that corresponds to nonnegative
  arguments).
  Note that both inequalities hold without any additional conditions
  on $\lambda$ or $\delta$.
  We only need to clarify it for the lower bound.
  Indeed, if $M \lambda < \omega_*^{-1}(M^2 \delta)$, then, in particular, we
  have $M \lambda < 1$ (since $\omega_*^{-1}$ takes values in
  $\ClosedOpenInterval{0}{1}$), and hence we can apply the bound
  $M^2 \delta \leq \omega_*(M \lambda) < M^2 \delta$, which leads to a
  contradiction.

  Applying \cref{St:Omega-Bounds}, we obtain that
  \[
    \omega^{-1}(M^2 \delta)
    \leq
    \sqrt{2 M^2 \delta} + M^2 \delta
    =
    M \sqrt{2 \delta} + M^2 \delta
    =
    M \sqrt{2 \delta} (1 + \tfrac{1}{\sqrt{2}} M \sqrt{\delta} \,).
  \]
  Thus,
  \[
    \lambda
    \leq
    \sqrt{2 \delta} (1 + \tfrac{1}{\sqrt{2}} M \sqrt{\delta} \,)
    =
    \sqrt{2 \delta} (1 + \tfrac{1}{4} H \sqrt{\delta} \,).
  \]
  Squaring both sides, we obtain the upper bound in
  \cref{LocalMetrics:GradNorm-Via-Func}.

  On the other hand, applying \cref{St:Omega-Bounds}, we get
  \[
    \omega_*(M^2 \delta)
    \geq
    1 - e^{-a}
    =
    1 - \frac{1}{e^a}
    \geq
    1 - \frac{1}{1 + a}
    =
    \frac{a}{1 + a},
  \]
  where
  \[
    a
    \DefinedEqual
    \sqrt{2 M^2 \delta} + \tfrac{1}{3} M^2 \delta
    =
    M \sqrt{2 \delta} + \tfrac{1}{3} M^2 \delta
    =
    M \sqrt{2 \delta} (1 + \tfrac{1}{3 \sqrt{2}} M \sqrt{\delta} \,).
  \]
  Therefore,
  \[
    \omega_*(M^2 \delta)
    \geq
    \frac{
      M \sqrt{2 \delta} (1 + \frac{1}{3 \sqrt{2}} M \sqrt{\delta} \,)
    }{
      1 + M \sqrt{2 \delta} (1 + \frac{1}{3 \sqrt{2}} M \sqrt{\delta} \,)
    }
    \geq
    \frac{M \sqrt{2 \delta}}{1 + M \sqrt{2 \delta}}.
  \]
  Thus,
  \[
    \lambda
    \geq
    \sqrt{2 \delta} (1 + M \sqrt{2 \delta} \,)
    =
    \sqrt{2 \delta} (1 + \tfrac{1}{2} H \sqrt{\delta} \,).
  \]
  Squaring both sides, we obtain the lower bound in
  \cref{LocalMetrics:GradNorm-Via-Func}.

  \ProofPart

  Now let us prove \cref{LocalMetrics:GradNormOpt-Via-Func}.
  We have
  \begin{gather*}
    2 \delta = \DualPairing{J_2 h}{h},
    \qquad
    h \DefinedEqual x - x_*,
    \qquad
    J_2 \DefinedEqual 2 \int_0^1 (1 - t) \Hessian f(x_* + t h) d t,
    \\
    \lambda_*^2
    =
    \DualPairing{\Gradient f(x)}{[\Hessian f(x_*)]^{-1} \Gradient f(x)}
    =
    \DualPairing{J_1 [\Hessian f(x_*)]^{-1} J_1 h}{h},
    \quad
    J_1 \DefinedEqual \int_0^1 \Hessian f(x^* + t h) d t.
  \end{gather*}
  Thus, we need to compare the operators $J_1 [\Hessian f(x_*)]^{-1} J_1$ and
  $J_2$.
  We do this as follows.
  First, we compare $J_1$ and $\Hessian f(x_*)$ using
  \cref{St:AuxiliaryBounds-Hess:Int}:
  \[
    (1 + \tfrac{1}{2} H \sqrt{\delta} \,)^{-1} \Hessian f(x_*)
    \preceq
    J_1
    \preceq
    (1 + \tfrac{1}{2} H \sqrt{\delta} \,) \Hessian f(x_*).
  \]
  From this, it follows that
  \[
    (1 + \tfrac{1}{2} H \sqrt{\delta} \,)^{-1} J_1
    \preceq
    J_1 [\Hessian f(x_*)]^{-1} J_1
    \preceq
    (1 + \tfrac{1}{2} H \sqrt{\delta} \,) J_1.
  \]
  Now it remains to compare $J_1$ and $J_2$:
  \begin{equation}
    \label{Proof:LocalMetrics:IntHess-Comp12}
    (1 + \tfrac{1}{6} H \sqrt{\delta} \,)^{-1} J_1
    \preceq
    J_2
    \preceq
    (1 + \tfrac{1}{6} H \sqrt{\delta} \,) J_1.
  \end{equation}
  Combining the above two relations, we obtain
  \cref{LocalMetrics:GradNormOpt-Via-Func}.

  To justify \cref{Proof:LocalMetrics:IntHess-Comp12}, first, note that
  \[
    \begin{aligned}
      J_2 - J_1
      &=
      \int_0^1 (1 - 2 t) \Hessian f(x^* + t h) d t
      \\
      &=
      \int_0^{1 / 2} (1 - 2 t) \Hessian f(x^* + t h) d t
      -
      \int_{1 / 2}^1 (2 t - 1) \Hessian f(x^* + t h) d t
      \\
      &=
      \int_0^{1 / 2} (1 - 2 t) \bigl[
        \Hessian f(x^* + t h)
        -
        \Hessian f\bigl( x^* + (1 - t) h \bigr)
      \bigr] d t,
    \end{aligned}
  \]
  where the last identity follows from the change of variables $t' = 1 - t$ in
  the second integral.
  Hence, by the Lipschitz continuity of the Hessian,
  \[
    \begin{aligned}
      \Norm{J_2 - J_1}
      &\leq
      \int_0^{1 / 2} (1 - 2 t) \Norm{
        \Hessian f(x^* + t h)
        -
        \Hessian f\bigl( x^* + (1 - t) h \bigr)
      } d t
      \\
      &\leq
      L_2 \Norm{h} \int_0^{1 / 2} (1 - 2 t)^2 d t
      =
      \frac{L_2}{2} \Norm{h} \int_0^1 t^2 d t
      =
      \frac{L_2}{6} \Norm{h},
    \end{aligned}
  \]
  and \cref{Proof:LocalMetrics:IntHess-Comp12} follows by strong convexity.
\end{proof}

Let us denote
\[
  \psi_0^L
  \DefinedEqual
  \psi\bigl( \Hessian f(x_*), G_0 \bigr).
\]
Also,
\[
  \sigma_k
  \DefinedEqual
  \sqrt{f(x_k) - f^*}.
\]
In what follows, we assume that, for all $k \geq 1$, we have
\begin{equation}
  \label{Bounded-Sum-Residuals}
  S_k \DefinedEqual H \sum_{i = 0}^{k - 1} \sigma_i \leq S.
\end{equation}

\begin{lemma}
  \label{St:Sum-Phi-Opt-Ubd}
  In \cref{Alg:BFGS}, for all $k \geq 1$, we have
  \[
    \sum_{i = 0}^{k - 1} \phi_{x_*}\bigl( d_i, \Gradient f(x_i) \bigr)
    \leq
    \psi_0^L + S_k + S_{k + 1}.
  \]
\end{lemma}

\begin{proof}
  Applying \cref{St:BregmanDist-Upd} with $C \DefinedEqual \Hessian f(x_*)$, we
  obtain, for all $i \geq 0$,
  \[
    \psi\bigl( \Hessian f(x_*), G_i \bigr)
    -
    \psi\bigl( \Hessian f(x_*), G_{i + 1} \bigr)
    =
    \phi_{x_*}(\delta_i, G_i \delta_i) - \phi_{x_*}(\delta_i, \gamma_i).
  \]
  Let $k \geq 1$ be arbitrary.
  Summing up the above inequalities for all $i = 0, \dots, k - 1$ and using the
  fact that $\psi(\cdot, \cdot) \geq 0$, we obtain
  \begin{equation}
    \label{Proof:Sum-Phi-Opt-Ubd:Prel}
    \sum_{i = 0}^{k - 1} \phi_{x_*}(\delta_i, G_i \delta_i)
    \leq
    \psi_0^L + \sum_{i = 0}^{k - 1} \phi_{x_*}(\delta_i, \gamma_i).
  \end{equation}
  Let us now estimate the sum in the right-hand side from above.

  Let $0 \leq i \leq k - 1$ be an arbitrary index.
  Note that
  \[
    \gamma_i = J_i \delta_i,
  \]
  where
  $
    J_i
    \DefinedEqual
    \int_0^1 \Hessian f\bigl( x_i + t (x_{i + 1} - x_i) \bigr) dt
  $.
  By \cref{St:AuxiliaryBounds-Hess:Int}, we have
  \[
    (1 + H \bar{\sigma}_i)^{-1} \Hessian f(x_*)
    \preceq
    J_i
    \preceq
    (1 + H \bar{\sigma}_i) \Hessian f(x_*),
  \]
  where
  \[
    \bar{\sigma}_i
    \DefinedEqual
    \frac{1}{2} (\sigma_i + \sigma_{i + 1}).
  \]
  Hence, in view of \cref{Def-Alpha-Beta},
  \begin{gather*}
    \alpha_{x_*}(\delta_i, \gamma_i)
    =
    \frac{
      \RelativeDualNorm{\gamma_i}{x_*}^2
    }{
      \DualPairing{\gamma_i}{\delta_i}
    }
    =
    % \begin{noindent}
    \frac{
      \DualPairing{J_i [\Hessian f(x_*)]^{-1} J_i \delta_i}{\delta_i}
    }{
      \DualPairing{J_i \delta_i}{\delta_i}
    }
    % \end{noindent}
    \leq
    1 + H \bar{\sigma}_i,
    \\
    \beta_{x_*}(\delta_i, \gamma_i)
    =
    \frac{
      \DualPairing{\gamma_i}{\delta_i}
    }{
      \RelativeNorm{\delta_i}{x_*}^2
    }
    =
    \frac{
      \DualPairing{J_i \delta_i}{\delta_i}
    }{
      \DualPairing{\Hessian f(x_*) \delta_i}{\delta_i}
    }
    \geq
    (1 + H \bar{\sigma}_i)^{-1}.
  \end{gather*}
  Therefore, according to \cref{Def-phi,Def-phi-Short},
  \[
    \phi_{x_*}(\delta_i, \gamma_i)
    =
    \alpha_{x_*}(\delta_i, \gamma_i) - 1
    -
    \ln \beta_{x_*}(\delta_i, \gamma_i)
    \leq
    H \bar{\sigma}_i + \ln(1 + H \bar{\sigma}_i)
    \leq
    2 H \bar{\sigma}_i,
  \]
  where we have used the inequality $\ln(1 + t) \leq t$, valid for any $t > -1$.

  Thus,
  \[
    \sum_{i = 0}^{k - 1} \phi_{x_*}(\delta_i, \gamma_i)
    \leq
    2 H \sum_{i = 0}^{k - 1} \bar{\sigma}_i
    =
    H \sum_{i = 0}^{k - 1} (\sigma_i + \sigma_{i + 1})
    \leq
    S_k + S_{k + 1}.
  \]
  Substituting this bound into \cref{Proof:Sum-Phi-Opt-Ubd:Prel}, we obtain
  \[
    \sum_{i = 0}^{k - 1} \phi_{x_*}(\delta_i, G_i \delta_i)
    \leq
    \psi_0^L + S_k + S_{k + 1}.
  \]
  It remains to use that
  \[
    \phi_{x_*}(\delta_i, G_i \delta_i)
    =
    \phi_{x_*}(d_i, G_i d_i)
    =
    \phi_{x_*}\bigl( d_i, \Gradient f(x_i) \bigr)
  \]
  in view of \cref{Alpha-Beta-Phi-Homogeneity} and the fact that
  $\delta_i = -h_i d_i$, $d_i = G_i^{-1} \Gradient f(x_i)$
  (see \cref{Alg:BFGS}).
\end{proof}

From \cref{St:Sum-Phi-Opt-Ubd}, using assumption~\eqref{Bounded-Sum-Residuals},
we can obtain the convergence rate (in average) for the
\emph{Dennis--More ratio}:
\begin{equation}
  \label{Def-DM}
  \xi_k
  \DefinedEqual
  \frac{
    \RelativeDualNorm{[G_k - \Hessian f(x_*)] d_k}{x_*}
  }{
    \RelativeNorm{d_k}{x_*}
  }
  =
  \frac{
    \RelativeDualNorm{\Gradient f(x_k) - \Hessian f(x_*) d_k}{x_*}
  }{
    \RelativeNorm{d_k}{x_*}
  }.
\end{equation}
The key characteristic, influencing the rate of convergence, is
\begin{equation}
  \label{Def-Omega}
  \Omega
  \DefinedEqual
  \psi_0^L + 2 S.
\end{equation}
In what follows, by
$\Map{\omega^{-1}}{\ClosedOpenInterval{0}{+\infty}}{\RealField}$, we denote the
inverse function of the right branch (corresponding to nonnegative arguments)
of the function $\omega$.

\begin{lemma}
  \label{St:DM-Rate}
  For all $k_0 \geq 0$ and all $k \geq 1$, we have
  \begin{equation}
    \label{DM-Rate}
    \frac{1}{k} \sum_{i = k_0}^{k_0 + k - 1} \xi_i
    \leq
    \omega^{-1}\Bigl( \frac{\Omega}{k} \Bigr)
    \leq
    \sqrt{\frac{2 \Omega}{k}} + \frac{\Omega}{k}.
  \end{equation}
\end{lemma}

\begin{proof}
  By \cref{St:Sum-Phi-Opt-Ubd}, assumption~\eqref{Bounded-Sum-Residuals} and the
  definition of~$\Omega$, we have
  \[
    \sum_{i = 0}^{k - 1} \phi_{x_*}\bigl( d_i, \Gradient f(x_i) \bigr)
    \leq
    \psi_0^L + 2 S
    =
    \Omega.
  \]
  Recall, from \cref{Def-phi,Def-phi-Short}, that
  \[
    \phi_{x_*}\bigl( d_i, \Gradient f(x_i) \bigr)
    =
    \alpha_i - 1
    -
    \ln \beta_i,
  \]
  where, according to \cref{Def-Alpha-Beta},
  \begin{gather*}
    \alpha_i
    \DefinedEqual
    \alpha_{x_*}\bigl( d_i, \Gradient f(x_i) \bigr)
    =
    \frac{
      \RelativeDualNorm{\Gradient f(x_i)}{x_*}^2
    }{
      \DualPairing{\Gradient f(x_i)}{d_i}
    },
    \\
    \beta_i
    \DefinedEqual
    \beta_{x_*}\bigl( d_i, \Gradient f(x_i) \bigr)
    =
    \frac{
      \DualPairing{\Gradient f(x_i)}{d_i}
    }{
      \RelativeNorm{d_i}{x_*}^2
    }.
  \end{gather*}
  Recall from \cref{Alpha-Beta-Relation} that $\alpha_i \geq \beta_i$.
  It can be shown that, for any $\alpha \geq \beta > 0$, the following
  inequality holds (see Lemma~2.3
  in~\cite{Rodomanov.Nesterov-RatesSuperlinear-22}):
  \[
    \alpha - 1 - \ln \beta
    \geq
    \omega( \sqrt{\alpha \beta - 2 \beta + 1} \,).
  \]
  Note that
  \[
    \begin{aligned}
      \alpha_i \beta_i - 2 \beta_i + 1
      &=
      \frac{
        \RelativeDualNorm{\Gradient f(x_i)}{x_*}^2
      }{
        \RelativeNorm{d_i}{x_*}^2
      }
      -
      2 \frac{
        \DualPairing{\Gradient f(x_i)}{d_i}
      }{
        \RelativeNorm{d_i}{x_*}^2
      }
      +
      1
      \\
      &=
      \frac{
        \RelativeDualNorm{\Gradient f(x_i) - \Hessian f(x_*) d_i}{x_*}^2
      }{
        \RelativeNorm{d_i}{x_*}^2
      }
      =
      \xi_i^2.
    \end{aligned}
  \]
  Thus, we obtain
  \[
    \phi_{x_*}\bigl( d_i, \Gradient f(x_i) \bigr)
    \geq
    \omega(\xi_i).
  \]
  Therefore,
  \[
    \sum_{i = 0}^{k - 1} \omega(\xi_i)
    \leq
    \Omega.
  \]
  In particular,
  \[
    \sum_{i = k_0}^{k_0 + k - 1} \omega(\xi_i) \leq \Omega.
  \]
  Since $\omega$ is convex, it follows that
  \[
    \omega\Bigl( \frac{1}{k} \sum_{i = k_0}^{k_0 + k - 1} \xi_i \Bigr)
    \leq
    \frac{1}{k} \sum_{i = k_0}^{k_0 + k - 1} \omega(\xi_i)
    \leq
    \frac{\Omega}{k}.
  \]
  This proves the first inequality in \cref{DM-Rate}.
  The second one follows from \cref{St:Omega-Bounds}.
\end{proof}

We will measure the rate of local superlinear convergence using the following
quantities:
\[
  r_k \DefinedEqual \RelativeNorm{x_k - x_*}{x_*}.
\]

\begin{lemma}
  \label{St:Residual-Via-DM}
  In \cref{Alg:BFGS}, for all $k \geq 0$, we have
  \begin{equation}
    \label{Residual-Via-DM}
    r_{k + 1}
    \leq
    \xi_k (r_k + r_{k + 1})
    +
    (1 + \tfrac{1}{2} H \sigma_k) \Abs{1 - h_k} r_k
    +
    \tfrac{1}{2} H \sigma_k r_k.
  \end{equation}
\end{lemma}

\begin{remark}
  The second term can be improved up to
  \[
    \max\Set{
      (1 + \tfrac{1}{2} H \sigma_k) h_k - 1,
      1 - (1 + \tfrac{1}{2} H \sigma_k)^{-1} h_k
    }.
  \]
\end{remark}

\begin{proof}
  Let $k \geq 0$ be arbitrary.
  By the definition of $x_{k + 1}$ in \cref{Alg:BFGS} and the triangle
  inequality, we have
  \begin{equation}
    \label{Proof:Residual-Via-DM:Prel}
    \begin{aligned}
      r_{k + 1}
      &=
      \RelativeNorm{x_{k + 1} - x_*}{x_*}
      =
      \RelativeNorm{x_k - x_* - h_k d_k}{x_*}
      \\
      &\leq
      \RelativeNorm{
      x_k - x_* - h_k [\Hessian f(x_*)]^{-1} \Gradient f(x_k)
      }{x_*}
      +
      h_k \RelativeNorm{[\Hessian f(x_*)]^{-1} \Gradient f(x_k) - d_k}{x_*}.
    \end{aligned}
  \end{equation}
  Since $\RelativeNorm{h}{A} = \RelativeDualNorm{A h}{A}$
  for any $h \in \VectorSpace{E}$
  and any $A \in \PDLin(\VectorSpace{E}, \VectorSpace{E}\Dual)$,
  we have
  \begin{equation}
    \label{Proof:Residual-Via-DM:MainTerm}
    \RelativeNorm{[\Hessian f(x_*)]^{-1} \Gradient f(x_k) - d_k}{x_*}
    =
    \RelativeDualNorm{\Gradient f(x_k) - \Hessian f(x_*) d_k}{x_*}
    =
    \xi_k \RelativeNorm{d_k}{x_*}.
  \end{equation}
  Similarly,
  \begin{equation}
    \label{Proof:Residual-Via-DM:AuxTerm}
    \begin{aligned}
      \hspace{2em}&\hspace{-2em}
      \RelativeNorm{
      x_k - x_* - h_k [\Hessian f(x_*)]^{-1} \Gradient f(x_k)
      }{x_*}
      =
      \RelativeDualNorm{
        \Hessian f(x_*) (x_k - x_*) - h_k \Gradient f(x_k)
      }{x_*}
      \\
      &\leq
      \Abs{1 - h_k}
      \RelativeDualNorm{\Gradient f(x_k)}{x_*}
      +
      \RelativeDualNorm{
        \Hessian f(x_*) (x_k - x_*) - \Gradient f(x_k)
      }{x_*}.
    \end{aligned}
  \end{equation}
  Note that
  \[
    f(x_k) = W_k (x_k - x_*),
    \qquad
    W_k \DefinedEqual \int_0^1 \Hessian f\bigl( x_* + t (x_k - x_*) \bigr).
  \]
  By \cref{St:AuxiliaryBounds-Hess:Int}, we have
  \[
    % this comment makes formatter work correctly here %
    [(1 + \tfrac{1}{2} H \sigma_k)]^{-1} \Hessian f(x_*)
    \preceq
    W_k
    \preceq
    (1 + \tfrac{1}{2} H \sigma_k) \Hessian f(x_*).
  \]
  Hence,
  \begin{gather*}
    \RelativeDualNorm{\Gradient f(x_k)}{x_*}
    =
    \RelativeDualNorm{W_k (x_k - x_*)}{x_*}
    \leq
    (1 + \tfrac{1}{2} H \sigma_k) r_k,
    \\
    \RelativeDualNorm{
      \Hessian f(x_*) (x_k - x_*) - \Gradient f(x_k)
    }{x_*}
    =
    \RelativeDualNorm{(\Hessian f(x_*) - W_k) (x_k - x_*)}{x_*}
    \leq
    \tfrac{1}{2} H \sigma_k r_k.
  \end{gather*}
  Substituting these estimates into
  \cref{Proof:Residual-Via-DM:AuxTerm}, we obtain
  \[
    \RelativeNorm{
    x_k - x_* - h_k [\Hessian f(x_*)]^{-1} \Gradient f(x_k)
    }{x_*}
    \leq
    (1 + \tfrac{1}{2} H \sigma_k) \Abs{1 - h_k} r_k
    +
    \tfrac{1}{2} H \sigma_k r_k.
  \]
  Combining this bound and \cref{Proof:Residual-Via-DM:MainTerm} into
  \cref{Proof:Residual-Via-DM:Prel}, we get
  \[
    r_{k + 1}
    \leq
    h_k \RelativeNorm{d_k}{x_*} \xi_k
    +
    (1 + \tfrac{1}{2} H \sigma_k) \Abs{1 - h_k} r_k
    +
    \tfrac{1}{2} H \sigma_k r_k.
  \]
  It remains to note that
  \[
    h_k \RelativeNorm{d_k}{x_*}
    =
    \RelativeNorm{x_{k + 1} - x_k}{x_*}
    \leq
    r_k + r_{k + 1}
  \]
  by the definition of $x_{k + 1}$ in \cref{Alg:BFGS} and the triangle
  inequality.
\end{proof}

Using the monotonicity of $\sigma_k$, we can show that $r_{k + 1}$ is not too
big compared to $r_k$ when $\sigma_k$ is sufficiently small.

\begin{lemma}
  \label{St:NextResidual-NotTooBig}
  For all $k \geq 0$, we have
  \begin{equation}
    \label{NextResidual-NotTooBig}
    r_{k + 1}
    \leq
    (1 + \tfrac{1}{3} H \sigma_k) r_k.
  \end{equation}
\end{lemma}

\begin{proof}
  Applying \cref{St:LocalMetrics}, we obtain
  \begin{gather*}
    r_{k + 1}^2
    \leq
    (1 + \tfrac{1}{3} H \sigma_{k + 1}) 2 \sigma_{k + 1}^2,
    \\
    r_k^2 \geq (1 + \tfrac{1}{3} H \sigma_k)^{-1} 2 \sigma_k^2.
  \end{gather*}
  Using that $\sigma_{k + 1} \leq \sigma_k$ and combining the first inequality
  with the second one, we obtain
  \[
    r_{k + 1}^2
    \leq
    (1 + \tfrac{1}{3} H \sigma_k) 2 \sigma_k^2
    \leq
    (1 + \tfrac{1}{3} H \sigma_k)^2 r_k^2.
  \]
  Taking the square root, we obtain \cref{NextResidual-NotTooBig}.
\end{proof}

\Cref{St:NextResidual-NotTooBig} allows us to eliminate $r_{k + 1}$ from the
right-hand side of \cref{Residual-Via-DM}.

\begin{corollary}
  \label{St:Residual-Via-DM:Simp}
  For all $k \geq 0$, we have
  \[
    r_{k + 1}
    \leq
    \bigl[
      (2 + \tfrac{1}{3} H \sigma_k) \xi_k
      +
      (1 + \tfrac{1}{2} H \sigma_k) \Abs{1 - h_k}
      +
      \tfrac{1}{2} H \sigma_k
    \bigr] r_k.
  \]
\end{corollary}

We already know that $\sigma_k \to 0$ (\cref{Bounded-Sum-Residuals}) and
$\xi_k \to 0$ (\cref{St:DM-Rate}).
Thus, to prove the superlinear convergence of $r_k$, it remains to show that
the step sizes $h_k \to 1$.
More precisely, we need to derive some convergence rate of $\Abs{1 - h_k}$
(on average).

Let us first consider the backtracking line search strategy.

We start with establishing sufficient conditions for accepting the unit step
size during line search.
Recall our main quantities:
\begin{gather*}
  \alpha_k
  \DefinedEqual
  \alpha_*\bigl( d_i, \Gradient f(x_k) \bigr)
  =
  \frac{
    \RelativeDualNorm{\Gradient f(x_k)}{x_*}^2
  }{
    \DualPairing{\Gradient f(x_k)}{d_k}
  },
  \\
  \beta_k
  \DefinedEqual
  \beta_*\bigl( d_k, \Gradient f(x_k) \bigr)
  =
  \frac{
    \DualPairing{\Gradient f(x_k)}{d_k}
  }{
    \RelativeNorm{d_k}{x_*}^2
  }.
\end{gather*}
Recall that $\alpha_k \geq \beta_k > 0$.

The following result shows that the unit step size satisfies Armijo's condition
as long as $\beta_k$ is sufficiently large compared to $H \sigma_k$.

\begin{lemma}
  \label{St:Armijo-UnitStepAccepted}
  Let $H \sigma_k \leq \Delta$, $\beta_k \geq \beta(\Delta)$
  for some $k \geq 0$, where
  \begin{gather}
    \label{Def-beta-lbd}
    \beta(\Delta)
    \DefinedEqual
    \frac{1}{4 (1 - \eta_1)}
    \bigl(
      \sqrt{(1 + \Delta)^2 + 8 (1 - \eta_1) \tilde{\Delta}}
      +
      (1 + \Delta)
    \bigr),
    \\
    \tilde{\Delta}
    \DefinedEqual
    \tilde{\Delta}(\Delta)
    \DefinedEqual
    \tfrac{1}{3} \Delta
    (1 + \tfrac{1}{2} \Delta)^{1 / 2}
    (1 + \tfrac{1}{6} \Delta)^{1 / 2}.
  \end{gather}
  Then,
  \begin{equation}
    \label{Armijo-UnitStepAccepted}
    f(x_k) - f(x_k - d_k)
    \geq
    \eta_1 \DualPairing{\Gradient f(x_k)}{d_k}.
  \end{equation}
\end{lemma}

\begin{proof}
  By \cref{LipHess-FuncBnd}, we have
  \[
    f(x_k - d_k) - f(x_k)
    \leq
    -\DualPairing{\Gradient f(x_k)}{d_k}
    +
    \frac{1}{2} \DualPairing{\Hessian f(x_k) d_k}{d_k}
    +
    \frac{L_2}{6} \Norm{d_k}^3.
  \]

  According to \cref{St:AuxiliaryBounds-Hess},
  $\Hessian f(x_k) \preceq (1 + H \sigma_k) \Hessian f(x_*)$.
  Therefore,
  \[
    \DualPairing{\Hessian f(x_k) d_k}{d_k}
    \leq
    (1 + H \sigma_k) \RelativeNorm{d_k}{x_*}^2
    =
    (1 + H \sigma_k) \beta_k^{-1}
    \DualPairing{\Gradient f(x_k)}{d_k}.
  \]

  Also, by strong convexity, $\mu B \preceq \Hessian f(x_*)$.
  Hence,
  \[
    L_2 \Norm{d_k}^3
    \leq
    \frac{L_2}{\mu^{3 / 2}} \RelativeNorm{d_k}{x_*}^3.
  \]
  By the definitions of $\alpha_k$ and $\beta_k$, we have
  \[
    \RelativeNorm{d_k}{x_*}^3
    =
    \RelativeNorm{d_k}{x_*}^2 \RelativeNorm{d_k}{x_*}
    =
    \bigl[
      \beta_k^{-1} \DualPairing{\Gradient f(x_k)}{d_k}
    \bigr]
    \bigl[
      (\alpha_k \beta_k)^{-1 / 2}
      \RelativeDualNorm{\Gradient f(x_k)}{x_*}
    \bigr].
  \]
  Since $\alpha_k \geq \beta_k$, we have
  $(\alpha_k \beta_k)^{-1 / 2} \leq \beta_k^{-1}$.
  Also, by \cref{St:LocalMetrics},
  \[
    \RelativeDualNorm{\Gradient f(x_k)}{x_*}
    \leq
    \sqrt{2} \sigma_k
    (1 + \tfrac{1}{6} H \sigma_k)^{1 / 2}
    (1 + \tfrac{1}{2} H \sigma_k)^{1 / 2}.
  \]
  Thus,
  \[
    \RelativeNorm{d_k}{x_*}^3
    \leq
    \beta_k^{-2} \sqrt{2} \sigma_k
    (1 + \tfrac{1}{6} H \sigma_k)^{1 / 2}
    (1 + \tfrac{1}{2} H \sigma_k)^{1 / 2}
    \DualPairing{\Gradient f(x_k)}{d_k},
  \]
  and hence, in view of \cref{Def-H} and the definition of
  $\tilde{\Delta}(\cdot)$,
  \[
    \begin{aligned}
      L_2 \Norm{d_k}^3
      &\leq
      \beta_k^{-2} H \sigma_k
      (1 + \tfrac{1}{6} H \sigma_k)^{1 / 2}
      (1 + \tfrac{1}{2} H \sigma_k)^{1 / 2}
      \DualPairing{\Gradient f(x_k)}{d_k}
      \\
      &=
      3 \beta_k^{-2} \tilde{\Delta}(H \sigma_k)
      \DualPairing{\Gradient f(x_k)}{d_k}.
    \end{aligned}
  \]

  Combining everything, we obtain
  \[
    \begin{aligned}
      \hspace{2em}&\hspace{-2em}
      \frac{1}{2} \DualPairing{\Hessian f(x_k) d_k}{d_k}
      +
      \frac{L_2}{6} \Norm{d_k}^3
      \\
      &\leq
      \frac{1}{2} \beta_k^{-1} (1 + H \sigma_k)
      \DualPairing{\Gradient f(x_k)}{d_k}
      +
      \frac{1}{2} \beta_k^{-2} \tilde{\Delta}(H \sigma_k)
      \DualPairing{\Gradient f(x_k)}{d_k}
      \\
      &=
      \frac{1}{2} \beta_k^{-1}
      \DualPairing{\Gradient f(x_k)}{d_k}
      \bigl(
        1 + H \sigma_k + \beta_k^{-1} \tilde{\Delta}(H \sigma_k)
      \bigr).
    \end{aligned}
  \]
  Thus,
  \[
    f(x_k) - f(x_k - d_k)
    \geq
    \DualPairing{\Gradient f(x_k)}{d_k}
    \bigl[
      1
      -
      \tfrac{1}{2} \beta_k^{-1}
      \bigl(
        1 + H \sigma_k + \beta_k^{-1} \tilde{\Delta}(H \sigma_k)
      \bigr)
    \bigr].
  \]

  Consequently, to ensure \cref{Armijo-UnitStepAccepted}, we need to ensure
  \[
    \beta_k^{-1}
    \bigl(
      1 + H \sigma_k + \beta_k^{-1} \tilde{\Delta}(H \sigma_k)
    \bigr)
    \leq
    2 (1 - \eta_1).
  \]
  Since $H \sigma_k \leq \Delta$ and $\tilde{\Delta}(\cdot)$ is increasing, it
  suffices to ensure that
  \[
    \beta_k^{-1} (1 + \Delta + \beta_k^{-1} \tilde{\Delta})
    \leq
    2 (1 - \eta_1).
  \]
  Equivalently,
  \[
    2 (1 - \eta_1) \beta_k^2
    \geq
    (1 + \Delta) \beta_k + \tilde{\Delta}.
  \]
  Denote $q \DefinedEqual 4 (1 - \eta_1)$.
  In this notation, our inequality becomes
  \[
    \tfrac{1}{2} q \beta_k^2
    \geq
    (1 + \Delta) \beta_k + \tilde{\Delta}
    \iff
    \beta_k^2
    \geq
    2 q^{-1} (1 + \Delta) \beta_k + 2 q^{-1} \tilde{\Delta}.
  \]
  Equivalently,
  \[
    \bigl( \beta_k - q^{-1} (1 + \Delta) \bigr)^2
    \geq
    [q^{-1} (1 + \Delta)]^2 + 2 q^{-1} \tilde{\Delta}.
  \]
  Thus, we finally obtain the solution
  \[
    \beta_k
    \geq
    \sqrt{[q^{-1} (1 + \Delta)]^2 + 2 q^{-1} \tilde{\Delta}}
    +
    q^{-1} (1 + \Delta)
    =
    \frac{1}{q}
    \bigl(
      \sqrt{(1 + \Delta)^2 + 2 q \tilde{\Delta}}
      +
      (1 + \Delta)
    \bigr).
  \]
  This is exactly $\beta(\Delta)$.
\end{proof}

By \cref{St:Sum-Phi-Opt-Ubd}, we know that $\beta_k \to 1$.
Thus, to ensure that the unit step size is eventually always accepted, we need
to make sure that $\beta(\Delta) < 1$ in \cref{St:Armijo-UnitStepAccepted}.
Observe that
\[
  \beta(0)
  =
  \frac{2}{4 (1 - \eta_1)}
  =
  \frac{1}{2 (1 - \eta_1)}
  <
  1
  \quad \iff \quad
  \eta_1 < \frac{1}{2}.
\]
Thus, from now on, we need to assume that
$\eta_1 \in \OpenOpenInterval{0}{\frac{1}{2}}$.
Since $\beta(\cdot)$ is monotone, we can therefore guarantee that
$\beta(\Delta) < 1$ for sufficiently small $\Delta$.
Let us now establish specific bounds.

\begin{lemma}
  \label{St:beta-SuffSmall}
  Let $\eta_1 \in \OpenOpenInterval{0}{\frac{1}{2}}$, and let $\beta(\cdot)$ be
  defined by \cref{Def-beta-lbd}.
  Let $\hat{\beta} \in \ClosedClosedInterval{\beta(0)}{1}$.
  Then,
  \[
    \beta(\Delta) \leq \hat{\beta}
    \iff
    \hat{\beta} \Delta + \tilde{\Delta}
    \leq
    \hat{\beta} \Bigl( \frac{\hat{\beta}}{\beta(0)} - 1 \Bigr).
  \]
\end{lemma}

\begin{proof}
  According to \cref{Def-beta-lbd}, we have $\beta(\Delta) \leq \hat{\beta}$ iff
  \[
    \sqrt{(1 + \Delta)^2 + 8 (1 - \eta_1) \tilde{\Delta}}
    +
    (1 + \Delta)
    \leq
    4 (1 - \eta_1) \hat{\beta}.
  \]
  Note that, for any $a, b \geq 0$ and $c > 0$, we have
  \[
    \sqrt{a^2 + b} + a \leq c
    \iff
    2 a c + b \leq c^2.
  \]
  Indeed, if $\sqrt{a^2 + b} + a \leq c$, then
  $a^2 + b \leq (c - a)^2 = c^2 - 2 a c + a^2$.
  If $2 a c + b \leq c^2$, then $a \leq c$ and $a^2 + b \leq (c - a)^2$, and
  hence $\sqrt{a^2 + b} \leq c - a$.

  Thus, $\beta(\Delta) \leq \hat{\beta}$ iff
  \[
    8 (1 - \eta_1) (1 + \Delta) \hat{\beta}
    +
    8 (1 - \eta_1) \tilde{\Delta}
    \leq
    16 (1 - \eta_1)^2 \hat{\beta}^2.
  \]
  After cancellations, this is equivalent to
  \[
    (1 + \Delta) \hat{\beta}
    +
    \tilde{\Delta}
    \leq
    2 (1 - \eta_1) \hat{\beta}^2
    \equiv
    \frac{\hat{\beta}^2}{\beta(0)}.
  \]
  Rearranging, we obtain the claim.
\end{proof}

In what follows, we will need to upper bound the (negative) quantity
$\beta(\Delta) - 1$, assuming that $\Delta$ is sufficiently small.
Let us show how to approximate the ``ideal'' value $\beta(0) - 1$ with ``high''
relative accuracy.

\begin{lemma}
  \label{St:Approx-beta-minus-one-Prel}
  Let $\eta_1 \in \OpenOpenInterval{0}{\frac{1}{2}}$ and
  let $q \in \OpenOpenInterval{0}{1}$.
  Denote
  \[
    \nu(q)
    \DefinedEqual
    \frac{
      6 q
    }{
      \sqrt{6 \tau q (1 - 2 \eta_1) + (3 + \tau)^2}
      +
      3 + \tau,
    }
    \qquad
    \tau
    \DefinedEqual
    \tau(q)
    \DefinedEqual
    \frac{2 (1 - \eta_1)}{1 + (1 - 2 \eta_1) q}.
  \]
  Then:
  \[
    \Delta \leq \nu(q) (1 - 2 \eta_1)
    \implies
    \beta(\Delta) - 1
    \leq
    (1 - q) [\beta(0) - 1]
    =
    -(1 - q) \frac{1 - 2 \eta_1}{2 (1 - \eta_1)}.
  \]
\end{lemma}

\begin{proof}
  We need to ensure that
  \[
    \beta(\Delta)
    \leq
    \hat{\beta}
    \DefinedEqual
    (1 - q) \beta(0) + q.
  \]
  Note that $\hat{\beta} \in \OpenOpenInterval{\beta(0)}{1}$.
  By \cref{St:beta-SuffSmall}, we need to ensure that
  \[
    \hat{\beta} \Delta + \tilde{\Delta}
    \leq
    \hat{\beta} \Bigl( \frac{\hat{\beta}}{\beta(0)} - 1 \Bigr)
    \EqualDefines
    \rho.
  \]
  Note that
  \[
    \tilde{\Delta}
    =
    \tfrac{1}{3} \Delta
    (1 + \tfrac{1}{2} \Delta)^{1 / 2}
    (1 + \tfrac{1}{6} \Delta)^{1 / 2}
    \leq
    \tfrac{1}{3} \Delta (1 + \tfrac{1}{2} \Delta)
    =
    \tfrac{1}{3} \Delta + \tfrac{1}{6} \Delta^2.
  \]
  Therefore,
  \[
    \hat{\beta} \Delta + \tilde{\Delta}
    \leq
    (\hat{\beta} + \tfrac{1}{3}) \Delta + \tfrac{1}{6} \Delta^2.
  \]
  Thus, it suffices to ensure that
  \[
    (\hat{\beta} + \tfrac{1}{3}) \Delta
    +
    \tfrac{1}{6} \Delta^2
    \leq
    \rho.
  \]
  This inequality is equivalent to
  \[
    \Delta^2 + 2 (3 \hat{\beta} + 1) \leq 6 \rho.
  \]
  Solving it, we obtain
  \[
    \Delta
    \leq
    \sqrt{6 \rho + (3 \hat{\beta} + 1)^2} - (3 \hat{\beta} + 1)
    =
    \frac{6 \rho}{\sqrt{6 \rho + (3 \hat{\beta} + 1)^2} + 3 \hat{\beta} + 1}
    \EqualDefines
    \bar{\Delta}.
  \]

  Note that
  \[
    \frac{\hat{\beta}}{\beta(0)} - 1
    =
    (1 - q) + \frac{q}{\beta(0)} - 1
    =
    q \Bigl( \frac{1}{\beta(0)} - 1 \Bigr)
    =
    q [2 (1 - \eta_1) - 1]
    =
    q (1 - 2 \eta_1).
  \]
  Hence, $\rho = \hat{\beta} q (1 - 2 \eta_1)$, and so
  \[
    \bar{\Delta}
    =
    \frac{
      6 \hat{\beta} q (1 - 2 \eta_1)
    }{
      \sqrt{
        6 \hat{\beta} q (1 - 2 \eta_1)
        +
        (3 \hat{\beta} + 1)^2
      }
      +
      3 \hat{\beta} + 1
    }
    =
    \nu (1 - 2 \eta_1),
  \]
  where
  \[
    \nu
    \DefinedEqual
    \frac{
      6 q
    }{
      \sqrt{
        6 \hat{\beta}^{-1} q (1 - 2 \eta_1)
        +
        (3 + \hat{\beta}^{-1})^2
      }
      +
      3 + \hat{\beta}^{-1}.
    }
  \]
  It remains to note that
  \begin{equation}
    \label{Proof:Approx-beta-minus-one:beta-hat-via-tau}
    \hat{\beta}
    =
    (1 - q) \frac{1}{2 (1 - \eta_1)} + q
    =
    \frac{1 - q + 2 (1 - \eta_1) q}{2 (1 - \eta_1)}
    =
    \frac{1 + (1 - 2 \eta_1) q}{2 (1 - \eta_1)}
    =
    \tau^{-1}.
  \end{equation}
\end{proof}

Note that the coefficient $\nu$ in \cref{St:Approx-beta-minus-one-Prel} is
monotonically increasing in $\eta_1$.
This follows from the fact that $\tau$ is decreasing in $\eta_1$
(see \cref{Proof:Approx-beta-minus-one:beta-hat-via-tau}).
Thus, we can lower bound $\nu(q)$ by substituting $\eta_1 = 0$:
\[
  \begin{aligned}
    \nu(q)
    &\geq
    % \begin{noindent}
    \frac{
      6 q
    }{
      \sqrt{6 [2 / (1 + q)] q + (3 + [2 / (1 + q)])^2}
      +
      3 + [2 / (1 + q)]
    }
    % \end{noindent}
    \\
    &=
    \frac{
      6 q (1 + q)
    }{
      \sqrt{12 q (1 + q) + (3 (1 + q) + 2)^2}
      +
      3 (1 + q) + 2
    }
    \\
    &=
    \frac{
      6 q (1 + q)
    }{
      \sqrt{12 q (1 + q) + (5 + 3 q)^2}
      +
      (5 + 3 q)
    }
    =
    \frac{
      6 (1 + q)
    }{
      \sqrt{25 + 21 q (2 + q)} + 5 + 3 q
    } q
    \geq
    \frac{3}{5} q,
  \end{aligned}
\]
where the final inequality follows from the fact that the fraction in front
of~$q$ is monotonically increasing.
Indeed, the inverse of this fraction is
\[
  \begin{aligned}
    \frac{\sqrt{25 + 21 q (2 + q)} + 5 + 3 q}{1 + q}
    &=
    \sqrt{\frac{25 + 21 q (2 + q)}{(1 + q)^2}}
    +
    \frac{5 + 3 q}{1 + q}
    \\
    &=
    \sqrt{21 + \frac{4}{(1 + q)^2}}
    +
    3 + \frac{2}{1 + q},
  \end{aligned}
\]
which is monotonically decreasing in $q$.

On the other hand, we can upper bound $\nu(q)$ by substituting
$\eta_1 = \frac{1}{2}$:
\[
  \nu(q)
  \leq
  \frac{6 q}{2 (3 + 1)}
  =
  \frac{3}{4} q.
\]
We thus conclude that there is not much reason to keep the ``complex''
coefficient $\nu(q)$ in \cref{St:Approx-beta-minus-one-Prel}.
Let us therefore present a simplified version of
\cref{St:Approx-beta-minus-one-Prel}.

\begin{lemma}
  \label{St:Approx-beta-minus-one}
  Let $\eta_1 \in \OpenOpenInterval{0}{\frac{1}{2}}$ and
  let $q \in \OpenOpenInterval{0}{1}$.
  Then:
  \[
    \Delta \leq \frac{3}{5} q (1 - 2 \eta_1)
    \implies
    \beta(\Delta) - 1
    \leq
    (1 - q) [\beta(0) - 1]
    =
    -(1 - q) \frac{1 - 2 \eta_1}{2 (1 - \eta_1)}.
  \]
  In particular, choosing $q = 0.1$, we obtain
  \[
    \Delta \leq 0.06 (1 - 2 \eta_1)
    \implies
    \beta(\Delta) - 1
    \leq
    -0.9 \frac{1 - 2 \eta_1}{2 (1 - \eta_1)}.
  \]
\end{lemma}

We are now ready to prove the rate of superlinear convergence for the
backtracking line search.

\begin{theorem}
  \label{St:SuperlinearRate:Backtrack-Unit}
  Consider \cref{Alg:BFGS} with the backtracking line search described in
  \cref{Backtrack-Unit:Quants,Backtrack-Unit:Test,Backtrack-Unit:FinalQuants},
  where $\eta_1 \in \OpenOpenInterval{0}{\frac{1}{2}}$ is a certain (absolute)
  constant.
  Let $K_0^L \geq 0$ be an integer for which the point $x_{K_0^L}$ is
  sufficiently good:
  \[
    H \sigma_{K_0^L}
    \leq
    \Delta
    \DefinedEqual
    0.06 (1 - 2 \eta_1).
  \]
  Then, for all $k \geq 1$, we have
  \begin{equation}
    \label{SuperlinearRate:Backtrack-Unit}
    \begin{aligned}
      r_{K_0^L + k}
      &\leq
      \Bigl[
        2.02 \, \omega^{-1}\Bigl( \frac{\Omega}{k} \Bigr)
        +
        \frac{1}{k} (\Upsilon \Omega + S)
      \Bigr]^k r_{K_0^L}
      \\
      &\leq
      \Bigl[
        2.02 \, \sqrt{\frac{2 \Omega}{k}}
        +
        \frac{1}{k} \bigl( (\Upsilon + 2.02) \Omega + S \bigr)
      \Bigr]^k r_{K_0^L},
    \end{aligned}
  \end{equation}
  where
  \begin{equation}
    \label{Def-Upsilon}
    \Upsilon
    \DefinedEqual
    [\omega_*(\bar{t})]^{-1}
    \ (\geq [\omega_*(0.45)]^{-1} \approx 6.76 \geq 6),
    \qquad
    \bar{t}
    \DefinedEqual
    0.9 \frac{1 - 2 \eta_1}{2 (1 - \eta_1)}
    \ (\leq 0.45).
  \end{equation}
\end{theorem}

\begin{remark}
  The term ``$+S$'' can be improved up to
  $
    \sum_{k = K_0^L}^\infty \sigma_k
    \sim
    H \sigma_0 \kappa^2
    \sim
    \Delta \kappa^2
  $.
\end{remark}

\begin{proof}
  By \cref{St:Residual-Via-DM:Simp}, for all $k \geq 0$, we have
  \[
    r_{k + 1} \leq \rho_k r_k,
  \]
  where
  \[
    \begin{aligned}
      \rho_k
      &\DefinedEqual
      (2 + \tfrac{1}{3} H \sigma_k) \xi_k
      +
      (1 + \tfrac{1}{2} H \sigma_k) (1 - h_k)
      +
      \tfrac{1}{2} H \sigma_k
      \\
      &\leq
      (2 + \tfrac{1}{3} H \sigma_k) \xi_k
      +
      (1 - h_k)
      +
      H \sigma_k.
    \end{aligned}
  \]
  (We have used the fact that $h_k \in \OpenClosedInterval{0}{1}$ for the
  backtracking line search.)
  Therefore, applying the arithmetic-geometric mean inequality,
  for all $k \geq 0$, we obtain
  \[
    r_{K_0^L + k}
    \leq
    r_{K_0^L} \prod_{i = K_0^L}^{K_0^L + k - 1} \rho_i
    \leq
    \Bigl[ \frac{1}{k} \sum_{i = K_0^L}^{K_0^L + k - 1} \rho_i \Bigr]^k
    r_{K_0^L}.
  \]
  Let us estimate the expression in the brackets.

  Since the method is monotone, for all $i \geq K_0^L$, we have
  \[
    H \sigma_i \leq H \sigma_{K_0^L} \leq \Delta.
  \]
  Therefore,
  \[
    \sum_{i = K_0^L}^{K_0^L + k - 1} \rho_i
    \leq
    (2 + \tfrac{1}{3} \Delta) \sum_{i = K_0^L}^{K_0^L + k - 1} \xi_i
    +
    \sum_{i = K_0^L}^{K_0^L + k - 1} (1 - h_i)
    +
    H \sum_{i = K_0^L}^{K_0^L + k - 1} \sigma_i.
  \]
  Let us now bound the three terms separately.

  For the third term, we can use \cref{Bounded-Sum-Residuals}:
  \[
    H \sum_{i = K_0^L}^{K_0^L + k - 1} \sigma_i \leq S.
  \]

  For the first term, we use \cref{St:DM-Rate}:
  \[
    \sum_{i = K_0^L}^{K_0^L + k - 1} \xi_i
    \leq
    k \omega^{-1}\Bigl( \frac{\Omega}{k} \Bigr).
  \]

  Finally, let us estimate the second term.
  For this, let us define the following index set:
  \[
    I_k
    \DefinedEqual
    \SetBuilder{
      K_0^L \leq i \leq K_0^L + k - 1
    }{
      \beta_i < \beta(\Delta)
    },
  \]
  where $\beta(\Delta)$ is defined in \cref{Def-beta-lbd}.
  Thus, for any $K_0^L \leq i \leq K_0^L + k - 1$, such that $i \notin I_k$,
  we have $H \sigma_i \leq \Delta$ (since $i \geq K_0^L$) and
  $\beta_i \geq \beta(\Delta)$ (since $i \notin I_k$).
  Therefore, by \cref{St:Armijo-UnitStepAccepted}, for any such $i$, the unit
  step size is accepted by the line search, i.e., $h_i = 1$.
  Consequently,
  \[
    \sum_{i = K_0^L}^{K_0^L + k - 1} (1 - h_i)
    =
    \sum_{i \in I_k} (1 - h_i)
    \leq
    \Cardinality{I_k},
  \]
  where $\Cardinality{I_k}$ is the number of elements in $I_k$.
  Let us show that $\Cardinality{I_k}$ is uniformly upper bounded.
  Applying \cref{St:Sum-Phi-Opt-Ubd} and the fact that $\alpha_i \geq \beta_i$,
  we obtain
  \[
    \Omega
    \geq
    \sum_{i = 0}^{K_0^L + k - 1} (\alpha_i - 1 - \ln \beta_i)
    \geq
    \sum_{i = 0}^{K_0^L + k - 1} (\beta_i - 1 - \ln \beta_i)
    =
    \sum_{i = 0}^{K_0^L + k - 1} \omega(\beta_i - 1)
    \geq
    \sum_{i \in I_k} \omega(\beta_i - 1).
  \]
  By the definition of $I_k$ and \cref{St:Approx-beta-minus-one}, for our choice
  of $\Delta$, we have, for any $i \in I_k$,
  \[
    \beta_i - 1
    \leq
    \beta(\Delta) - 1
    \leq
    -\bar{t}.
  \]
  Since $\omega(\cdot)$ is decreasing for negative arguments, we therefore have
  $\omega(\beta_i - 1) \geq \omega(-\bar{t}) = \omega_*(\bar{t})$.
  Thus, we obtain
  \[
    \Omega \geq \omega_*(\bar{t}) \Cardinality{I_k},
  \]
  and hence
  \[
    \Cardinality{I_k} \leq \Upsilon \Omega.
  \]
  Thus,
  \[
    \sum_{i = K_0^L}^{K_0^L + k - 1} (1 - h_i)
    \leq
    \Upsilon \Omega.
  \]

  Putting the three bounds together, we obtain
  \[
    \frac{1}{k} \sum_{i = K_0^L}^{K_0^L + k - 1} \rho_i
    \leq
    (2 + \tfrac{1}{3} \Delta)
    \omega^{-1}\Bigl( \frac{\Omega}{k} \Bigr)
    +
    \frac{1}{k} (\Upsilon \Omega + S).
  \]
  It remains to note that
  \[
    2 + \tfrac{1}{3} \Delta
    =
    2 + \tfrac{1}{3} 0.06 (1 - 2 \eta_1)
    =
    2 + 0.02 (1 - 2 \eta_1)
    \leq
    2.02.
  \]
  This proves the first inequality in \cref{SuperlinearRate:Backtrack-Unit}.
  The second one follows from \cref{St:Omega-Bounds}.
\end{proof}

Let us discuss the efficiency estimate from
\cref{St:SuperlinearRate:Backtrack-Unit}.
From \cref{Def-Omega}, we see that $S \leq \Omega$.
Therefore, bound \cref{SuperlinearRate:Backtrack-Unit} becomes
\[
  r_{K_0^L + k}
  \leq
  \Bigl[
    2.02 \omega^{-1}\Bigl( \frac{\Omega}{k} \Bigr)
    +
    \frac{(\Upsilon + 1) \Omega}{k} \,
  \Bigr]^k r_{K_0^L}
  \leq
  \Bigl[
    2.02 \sqrt{\frac{2 \Omega}{k}}
    +
    \frac{(\Upsilon + 3.02) \Omega}{k}
  \Bigr]^k r_{K_0^L}.
\]
Thus, the main quantity in this bound is $\Omega$.
(Recall that $\Upsilon$ is some absolute constant depending on the line search
parameter $\eta_1$.)

\begin{corollary}
  \label{St:StartMomentOfSuperConv}
  Denote
  \[
    K_1^L
    \DefinedEqual
    \Ceil{7 (\Upsilon + 3.02) \Omega}.
  \]
  Then, for all $k \geq K_1^L$, we have
  \[
    r_{K_0^L + k} \leq 2^{-k} r_{K_0^L}.
  \]
\end{corollary}

\begin{proof}
  Indeed, let $k \geq K_1^L$.
  Then
  \[
    \frac{(\Upsilon + 3.02) \Omega}{k}
    \leq
    \frac{1}{7}
    \approx
    0.14.
  \]
  On the other hand, by \cref{Def-Upsilon}, we have
  \[
    2.02 \sqrt{\frac{2 \Omega}{k}}
    \leq
    2.02 \sqrt{\frac{2}{7 (\Upsilon + 3.02)}}
    \leq
    2.02 \sqrt{\frac{2}{7 (6.7 + 3.02)}}
    \approx
    0.34.
  \]
  Thus,
  \[
    2.02 \sqrt{\frac{2 \Omega}{k}}
    +
    \frac{(\Upsilon + 3.02) \Omega}{k}
    \lessapprox
    0.14 + 0.34
    =
    0.48
    \leq
    \frac{1}{2}.\qedhere
  \]
\end{proof}

\begin{remark}
  Recall that $\Upsilon$ is the constant which depends on the line search
  parameter $\eta_1 \in \OpenOpenInterval{0}{\frac{1}{2}}$.

  One reasonable value is
  \[
    \eta_1 = \frac{1}{4}.
  \]
  In this case,
  \begin{gather*}
    \bar{t}
    =
    0.9 \frac{1 - 1/2}{2 (1 - 1/4)}
    =
    0.9 \frac{1}{3}
    =
    0.3,
    \\
    \Upsilon
    =
    [\omega_*(0.3)]^{-1}
    \approx
    17.64
    \leq
    18.
  \end{gather*}

  Note that we can slightly decrease $\Upsilon$ by decreasing $\eta_1$.
  For example, for
  \[
    \eta_2 = \frac{1}{8},
  \]
  we have
  \[
    \bar{t} \approx 0.38,
    \qquad
    \Upsilon \approx 9.84 \leq 10.
  \]
\end{remark}

\begin{remark}
  The constant ``$7$'' in \cref{St:StartMomentOfSuperConv} can be slightly
  improved by making it dependent on $\Upsilon$.
\end{remark}

Let us estimate from above the main quantity
\[
  \Omega = \psi_0^L + 2 S.
\]

We start with estimating $\psi_0^L$.
Suppose that the initial Hessian approximation is chosen in the following way:
\[
  G_0 \DefinedEqual \mu_0 B,
\]
where $\mu_0 \in \ClosedClosedInterval{\mu}{L}$ is a certain constant.
Recall that, by \cref{StrongConv,LipGrad}, we have
\[
  \mu B \preceq \Hessian f(x_*) \preceq L B.
\]
Therefore,
\[
  \frac{\mu}{\mu_0} G_0
  \preceq
  \Hessian f(x_*)
  \preceq
  \frac{L}{\mu_0} G_0.
\]
Thus, the relative eigenvalues $\lambda_1, \dots, \lambda_n$ of~$G_0$
w.r.t.~$\Hessian f(x_*)$ can be estimated as follows:
\[
  \frac{\mu_0}{L} \preceq \lambda_i \leq \frac{\mu_0}{\mu},
  \qquad
  1 \leq i \leq n.
\]
Thus, according to \cref{BregmanDist-ViaEigs},
\[
  \psi_0^L
  =
  \psi\bigl( \Hessian f(x_*), G_0 \bigr)
  =
  \sum_{i = 1}^n \omega(\lambda_i - 1)
  \leq
  n \max\Set[\Big]{
    \omega\Bigl( \frac{\mu_0}{L} - 1 \Bigr),
    \omega\Bigl( \frac{\mu_0}{\mu} - 1 \Bigr)
  }.
\]
Since $\mu_0 \in \ClosedClosedInterval{\mu}{L}$, we have
\[
  \kappa^{-1} - 1 \leq \frac{\mu_0}{L} - 1 \leq 0,
  \qquad
  0 \leq \frac{\mu_0}{\mu} - 1 \leq \kappa - 1.
\]
Since $\omega$ is decreasing for negative arguments and increasing for positive
ones,
\[
  \omega\Bigl( \frac{\mu_0}{L} - 1 \Bigr)
  \leq
  \omega(\kappa^{-1} - 1),
  \qquad
  \omega\Bigl( \frac{\mu_0}{\mu} - 1 \Bigr)
  \leq
  \omega(\kappa - 1).
\]
Therefore,
\[
  \psi_0^L
  \leq
  n \max\Set{\omega(\kappa^{-1} - 1), \omega(\kappa - 1)}
  =
  n \omega(\kappa - 1)
  \leq
  n \kappa.
\]
(we have used that $\omega(t^{-1} - 1) \leq \omega(t - 1)$ for any $t \geq 1$).

Now let us estimate $S$.
According to \cref{St:LinearConvergence-Backtrack-Unit}, for all $k \geq 1$,
we have
\[
  \sigma_k^2 \leq (1 - \eta_1 \kappa^{-2} \zeta_k)^k \sigma_0^2,
\]
where $\zeta_k \DefinedEqual \exp(-\psi_0^G / k)$
with $\psi_0^G \leq n \ln \kappa$
(see the discussion after \cref{St:LinearConvergence-Started}).
Thus, we can bound $S$ as follows.
Denote $K_0^G = \Ceil{2 \psi_0^G} \leq \Ceil{2 n \ln \kappa}$.
Then, for all $k \geq K_0^G$, we have
$\zeta_k \geq \exp(-1 / 2) \geq \frac{1}{2}$ and hence
\begin{equation}
  \label{LinearConv-FuncRes:Backtrack-Unit}
  \sigma_k^2 \leq (1 - \tfrac{1}{2} \eta_1 \kappa^{-2})^k \sigma_0^2.
\end{equation}
Therefore, for all $k \geq K_0^G$, we have
\[
  \sigma_k
  \leq
  (1 - \tfrac{1}{2} \eta_1 \kappa^{-2})^{k / 2} \sigma_0
  \leq
  (1 - \tfrac{1}{4} \eta_1 \kappa^{-2})^k \sigma_0,
\]
where we have used that $\sqrt{1 - \tau} \leq 1 - \frac{1}{2} \tau$.
On the other hand, for all $0 \leq k \leq K_0^G - 1$, we have
$\sigma_k \leq \sigma_0$.
Thus, according to \cref{Bounded-Sum-Residuals},
\[
  \begin{aligned}
    S
    &=
    H \sum_{i = 0}^\infty \sigma_i
    =
    H \sum_{i = 0}^{K_0^G - 1} \sigma_i
    +
    H \sum_{i = K_0^G}^\infty \sigma_i
    \leq
    K_0^G H \sigma_0
    +
    H \sigma_0
    \sum_{i = K_0^G}^\infty (1 - \tfrac{1}{4} \eta_1 \kappa^{-2})^i
    \\
    &\leq
    \Ceil{2 n \ln \kappa} H \sigma_0
    +
    4 \eta_1^{-1} \kappa^2 H \sigma_0
    =
    (4 \eta_1^{-1} \kappa^2 + \Ceil{2 n \ln \kappa}) H \sigma_0.
  \end{aligned}
\]

Thus, we finally obtain the following estimate:
\[
  \Omega
  \leq
  n \kappa
  +
  2 (4 \eta_1^{-1} \kappa^2 + \Ceil{2 n \ln \kappa}) H \sigma_0.
\]
Depending on how small $\sigma_0$ is, the principal term in the above bound can
either be the first one or the second one.

It is interesting that we can significantly improve the first term in the above
bound (corresponding to $\psi_0^L$) by restricting the choice of the
coefficient $\mu_0$ defining the initial Hessian approximation $G_0$.
Indeed, suppose that
\[
  \mu \leq \mu_0 \leq \nu \mu,
\]
where $\nu \in \ClosedClosedInterval{1}{\kappa}$ is a certain absolute
constant.
Then, we still have $\mu_0 \in \ClosedClosedInterval{\mu}{L}$, but now
\[
  \kappa^{-1} - 1 \leq \frac{\mu_0}{L} - 1 \leq 0,
  \qquad
  0 \leq \frac{\mu_0}{\mu} - 1 \leq \nu - 1.
\]
Therefore,
\[
  \omega\Bigl( \frac{\mu_0}{L} - 1 \Bigr)
  \leq
  \omega(\kappa^{-1} - 1),
  \qquad
  \omega\Bigl( \frac{\mu_0}{\mu} - 1 \Bigr)
  \leq
  \omega(\nu - 1),
\]
and hence
\[
  \begin{aligned}
    \psi_0^L
    &\leq
    n \max\Set{\omega(\kappa^{-1} - 1), \omega(\nu - 1)}
    \\
    &\leq
    n [\omega(\kappa^{-1} - 1) + \omega(\nu - 1)]
    \leq
    n [\ln \kappa + \omega(\nu - 1)].
  \end{aligned}
\]
Thus, we obtain
\[
  \Omega
  \leq
  n [\ln \kappa + \omega(\nu - 1)]
  +
  2 (4 \eta_1^{-1} \kappa^2 + \Ceil{2 n \ln \kappa}) H \sigma_0.
\]
Assuming that $\nu$ is an absolute constant (e.g., $\nu = 1$ or $\nu = 2$), the
first term in this expression is now $\BigO(n \ln \kappa)$.

\begin{remark}
  In particular, we have proved that the starting moment of superlinear
  convergence of BFGS with backtracking line search and the initial Hessian
  approximation $G_0 = \mu_0 B$ is $\BigO(n \ln \kappa)$ when
  $\mu \leq \mu_0 \leq \nu \mu$ for an absolute constant
  $\nu \in \ClosedClosedInterval{1}{\kappa}$, provided that the starting point
  is sufficiently good: $H \sigma_0 = \BigO(\kappa^{-2})$.
  This complements the result from~\cite{Rodomanov.Nesterov-NewResults-21} where
  it was shown that such a good starting moment of superlinear convergence can
  be obtained for BFGS by using the initial Hessian approximation $G_0 = L B$
  (and with unit step size at every iteration).
\end{remark}

  \section{Total Complexity Bound}

Let us now put together the global and local superlinear convergence
guarantees, which we have established before, and estimate the total complexity
of \cref{Alg:BFGS} for finding an $\epsilon$-approximate solution for
problem~\eqref{Prob}.

We consider the backtracking line search described in
\cref{Backtrack-Unit:Quants,Backtrack-Unit:Test,Backtrack-Unit:FinalQuants} and
suppose that
\[
  G_0 = \mu_0 B,
  \qquad
  \mu_0 \in \ClosedClosedInterval{\mu}{\nu \mu},
  \qquad
  \nu \in \ClosedClosedInterval{1}{\kappa}.
\]

From \cref{LinearConv-FuncRes:Backtrack-Unit}, we know that
\[
  \delta_k
  \DefinedEqual
  f(x_k) - f_*
  \leq
  (1 - \tfrac{1}{2} \eta_1 \kappa^{-2})^k \delta_0
  \leq
  \exp(-\tfrac{1}{2} \eta_1 \kappa^{-2} k) \delta_0
\]
for all $k \geq K_0^G$, where $K_0^G \leq 2 n \ln \kappa$.
Thus, to make $\delta_k \leq \epsilon$ for some $0 < \epsilon < \delta_0$, it
suffices to make $k \geq \max\Set{ K_0^G, K_G(\epsilon) }$ iterations, where
\[
  K_G(\epsilon)
  \DefinedEqual
  2 \eta_1^{-1} \kappa^2 \ln \frac{\delta_0}{\epsilon}.
\]
Thus, we can guarantee that $\delta_k \leq \epsilon$ after at most the
following number of iterations:
\[
  \Ceil{
    2 \eta_1^{-1} \kappa^2 \ln \frac{\delta_0}{\epsilon}
    +
    2 n \ln \kappa
  }.
\]

Let us now estimate the number of iterations $K_0^L$ needed to reach the region
of superlinear convergence as described in
\cref{St:SuperlinearRate:Backtrack-Unit}.
Thus, we need to make $H \sigma_k \leq \Delta$, or, equivalently,
$H^2 \delta_k \leq \Delta^2$.
For this, we need to take
$\epsilon = H^{-2} \Delta^2 = [0.06 (1 - 2 \eta_1)]^2 H^{-2}$.
This gives us
\[
  K_0^L
  =
  \Ceil[\Big]{
    2 \eta_1^{-1} \kappa^2
    \ln\bigl( [0.06 (1 - 2 \eta_1)]^{-2} H^2 \delta_0 \bigr)
    +
    2 n \ln \kappa
  }.
\]
Note that, for all $k \leq K_0^L$, we have
$H \sigma_k \leq H \sigma_{K_0}^L \leq \Delta$.

By \cref{St:StartMomentOfSuperConv}, for all $k \geq K_1^L$, where
\[
  K_1^L
  \DefinedEqual
  \Ceil{7 (\Upsilon + 3.02) \Omega},
\]
we have
\[
  r_{K_0^L + k} \leq 2^{-k} r_{K_0^L}.
\]
By \cref{St:LocalMetrics}, for all $k \geq 0$, we have
\[
  (1 + \tfrac{1}{3} H \sigma_k)^{-1} 2 \delta_k
  \leq
  r_k^2
  \leq
  (1 + \tfrac{1}{3} H \sigma_k) 2 \delta_k.
\]
Since $H \sigma_k \leq \Delta$ for all $k \geq K_0^L$, we therefore have, for
all $k \geq K_0^L$:
\[
  (1 + \tfrac{1}{3} \Delta)^{-1} \, 2 \delta_k
  \leq
  r_k^2
  \leq
  (1 + \tfrac{1}{3} \Delta) 2 \delta_k.
\]
In particular, for all $k \geq K_0^L$,
\[
  r_k^2 \geq (1 + \tfrac{1}{3} \Delta)^{-1} \, 2 \delta_k.
\]
Also, since $H^2 \delta_{K_0^L} = (H \sigma_{K_0^L})^2 \leq \Delta$, we have
\[
  r_{K_0^L}^2
  \leq
  (1 + \tfrac{1}{3} \Delta) 2 \delta_{K_0^L}
  \leq
  (1 + \tfrac{1}{3} \Delta) 2 H^{-2} \Delta^2.
\]
Thus, we obtain, for all $k \geq K_1^L$,
\[
  (1 + \tfrac{1}{3} \Delta)^{-1} \, 2 \delta_{K_0^L + k}
  \leq
  r_k^2
  \leq
  2^{-2 k} r_{K_0^L}^2
  \leq
  2^{-2 k} (1 + \tfrac{1}{3} \Delta) 2 H^{-2} \Delta^2.
\]
Consequently, for all $k \geq K_1^L$,
\[
  \delta_{K_0^L + k}
  \leq
  2^{-2 k} (1 + \tfrac{1}{3} \Delta)^2 H^{-2} \Delta^2.
\]
Thus, we can ensure that $\delta_{K_0^L + k} \leq \epsilon$
for any $k \geq \max\Set{K_1^L, K_1^L(\epsilon)}$, where
\[
  K_1^L(\epsilon)
  \DefinedEqual
  \frac{1}{2}
  \log_2 \frac{(1 + \tfrac{1}{3} \Delta)^2 \Delta^2}{H^2 \epsilon}.
\]
In our case, $\Delta = 0.06 (1 - 2 \eta_1) \leq 0.06$.
Therefore,
\[
  (1 + \tfrac{1}{3} \Delta)^2 \Delta^2
  \leq
  (1.02 \cdot 0.06)^2
  = 0.0037\ldots
  \leq
  0.01.
\]
Thus,
\[
  K_1^L(\epsilon)
  \leq
  \frac{1}{2} \log_2 \frac{0.01}{H^2 \epsilon}.
\]

To summarize: in order to make $\delta_k \leq \epsilon$, it suffices to make
the following number of iterations:
\[
  \begin{aligned}
    \hspace{2em}&\hspace{-2em}
    K_0^L + K_1^L + K_1^L(\epsilon)
    \\
    &\leq
    \Ceil[\Big]{
      2 \eta_1^{-1} \kappa^2
      \ln \frac{H^2 \delta_0}{[0.06 (1 - 2 \eta_1)]^2}
      +
      2 n \ln \kappa
    }
    +
    \Ceil{7 (\Upsilon + 3.02) \Omega}
    +
    \Ceil[\Big]{
      \frac{1}{2} \log_2 \frac{0.01}{H^2 \epsilon}
    }.
  \end{aligned}
\]
The first term is the number of iterations to reach the region of local
convergence.
The second one is the number of iterations to begin superlinear convergence.
The third one is the number of iterations to obtain $\epsilon$-approximate
solution at a superlinear rate.
  \section{BFGS with Restarts}

As can be seen from the complexity estimate, which we established in the
previous section, the standard BFGS Method has the following drawback.
Once it gets into the region of local convergence ($H \sigma_k \leq \Delta$),
it needs
\[
  K_1^L
  =
  \BigO(\Omega)
  =
  \BigO\bigl( \psi_0^L + (\kappa^2 + n \ln \kappa) H \sigma_0 \bigr),
\]
where $\BigO(\cdot)$ hides some absolute constants and the constants, depending
on the line search parameter $\eta_1$, and
\[
  \psi_0^L
  =
  \begin{cases}
    \BigO(n \kappa),
    & \text{if $\mu_0 \in \ClosedClosedInterval{\mu}{L}$},
    \\
    \BigO(n \ln \kappa + n),
    & \text{if $\mu_0 \in \ClosedClosedInterval{\mu}{\nu \mu}$},
  \end{cases}
\]
where $G_0 = \mu_0 B$ and $\nu \in \ClosedClosedInterval{1}{\kappa}$ is an
absolute constant.
We see that the second term in this estimate can be of the order
\[
  \BigO(\kappa^2 H \sigma_0),
\]
which is rather big when the initial functional residual
$f(x_0) - f_* \equiv \sigma_0^2$ is large.

Informally speaking, this happens because BFGS may ``spoil'' the initial
Hessian approximation $G_0 = \mu_0 B$ too much during the global convergence
phase, and then it arrives into the region of local convergence with the ``very
bad'' Hessian approximation $G_{K_0^L}$ which may be much worse than $G_0$.
Indeed, if the starting point~$x_0$ were already in the region of local
convergence, as described, e.g., by
\[
  (\kappa^2 + n \ln \kappa) H \sigma_0 = \BigO(1),
\]
then $K_1^L$ would be much smaller:
\[
  \tilde{K}_1^L \DefinedEqual \BigO(\psi_0^L).
\]

A natural approach to deal with this issue is of course to restart the method
with $G_0 = \mu_0 B$ when it arrives into the region of local convergence.
Then, we would have the total complexity of
\[
  \tilde{K}_0^L + \tilde{K}_1^L + K_1^L(\epsilon)
  =
  \BigO\bigl(
    \tilde{K}_0^L + \psi_0^L + \log(1 / [H^2 \epsilon])
  \bigr)
\]
iterations to find an $\epsilon$-approximate solution, where $\tilde{K}_0^L$ is
the number of iterations sufficient to ensure that
\[
  (\kappa^2 + n \ln \kappa) H \sigma_{\tilde{K}_0^L} = \BigO(1).
\]
It is not difficult to see that
\[
  \tilde{K}_0^L = \BigOTilde(\kappa^2 + n),
\]
where $\BigOTilde(\cdot)$ hides logarithmic factors depending on $\kappa$, $n$,
$H$, $\sigma_0$, etc.

However, direct implementation of this idea is not so simple.
First, we cannot measure the functional residual $f(x_k) - f_*$ since we may
not know the optimal value $f_*$.
Second, even if we could somehow efficiently estimate $f(x_k) - f_*$, we still
need to know other constants ($\kappa$, $H_2$, etc.) which enter the
description of the region of local convergence.
Finally, even if we also knew all these constants, our theoretical description
of the region of local convergence may be too pessimistic for real-world
problems.
Thus, ideally, we would like to use less rigid restarts than the direct one
outlined above.

Let us present one such procedure whose worst-case complexity is of the same
order as the complexity of the direct-restart approach.
The basic idea is very simple: given some number $N \geq 1$, we simply restart
the method after $N, 2 N, 4 N, 8 N, \ldots$ iterations.
A reasonable value of the parameter $N$ is $N = n$.

In what follows, for a point $y_0 \in \RealField^n$ and a number $N \geq 1$, by
$\BFGS(y_0, N)$, we denote the point, produced by the BFGS method after~$N$
iterations started from the initial point~$y_0$ (and the initial Hessian
approximation $G_0 = \mu_0 B$).

\begin{SimpleAlgorithm}[
    title = {BFGS with Restarts},
    label = Alg:BFGS-Restarts,
    width = 0.5\linewidth
  ]
  \begin{AlgorithmGroup}[Input]
    Initial point $y_0 \in \RealField^n$, number $N \geq 1$.
  \end{AlgorithmGroup}
  \AlgorithmGroupSeparator
  \begin{AlgorithmGroup}[Iteration $t \geq 0$]
    \begin{AlgorithmSteps}
      \AlgorithmStep
        Denote $N_t \DefinedEqual 2^t N$.
      \AlgorithmStep
        Compute $y_{t + 1} \DefinedEqual \BFGS(y_t, N_t)$.
    \end{AlgorithmSteps}
  \end{AlgorithmGroup}
\end{SimpleAlgorithm}

\begin{theorem}
  Consider \cref{Alg:BFGS-Restarts}.
  Let $K_0^L$, $\tilde{K}_1^L$ and $K_1^L(\epsilon)$ be defined as above but
  with $\sigma_0 \DefinedEqual \sqrt{f(y_0) - f_*}$.
  Then, the total number of inner-level iterations to obtain a point
  $f(y_t) - f_* \leq \epsilon$ is at most
  \[
    4 [ \tilde{K}_0^L + \tilde{K}_1^L + K_1^L(\epsilon) ]
    =
    \BigO\bigl(
      \tilde{K}_0^L + \psi_0^L + \log(1 / [H^2 \epsilon])
    \bigr).
  \]
  (Assuming that $N \leq n$.)
\end{theorem}

\begin{proof}
  \ProofPart

  Since the BFGS method is monotone, for all $t \geq 0$, we have
  \[
    f(y_{t + 1}) \leq f(y_t).
  \]

  \ProofPart

  Let $t_0 \geq 0$ be the smallest integer such that
  \[
    N_{t_0} \geq \tilde{K}_0^L.
  \]
  Consider the epoch with index~$t_0$. Since $\sigma(y_{t_0}) \leq \sigma_0$ and
  $N_{t_0} \geq K_0^L$, at the end of this epoch, we obtain the point
  $y_{t_0 + 1}$ in the region of local convergence:
  \[
    (\kappa^2 + n \ln \kappa) H \sigma(y_{t_0 + 1}) = \BigO(1).
  \]

  \ProofPart

  Let $t_1 \geq t_0 + 1$ be the smallest integer such that
  \[
    N_{t_1} \geq \tilde{K}_1^L + K_1^L(\epsilon),
  \]
  where
  \[
    \tilde{K}_1^L
    =
    \BigO\bigl(
      \psi_0^L + (\kappa^2 + n \ln \kappa) H \sigma(y_{t_0 + 1})
    \bigr)
    =
    \BigO(\psi_0^L).
  \]
  Clearly,
  \[
    N_{t_1}
    \leq
    2 [\tilde{K}_0^L + \tilde{K}_1^L + K_1^L(\epsilon)].
  \]

  Consider the epoch with index $t_1$.
  By the monotonicity, we have $\sigma(y_{t_1}) \leq \sigma(y_{t_0 + 1})$.
  Therefore, after at most $\tilde{K}_1^L$ iterations, the BFGS method will
  start converging at least with the fast linear rate $2^{-k}$.
  After that, it will need at most $K_1^L(\epsilon)$ iterations to find an
  $\epsilon$-approximate solution.
  Since $\tilde{K}_1^L + K_1^L(\epsilon) \leq N_{t_1}$, all of this will happen
  during the epoch~$t_1$ without any restarts.
  Since the inner method is monotone, at the end of epoch~$t_1$, we will thus
  have an $\epsilon$-approximate solution $y_{t_1 + 1}$:
  \[
    f(y_{t_1 + 1}) - f_* \leq \epsilon.
  \]

  \ProofPart

  It remains to estimate the total number of inner iterations.
  We need to do at most $t_1$ epochs.
  Therefore, the total number of inner iterations is
  \[
    \begin{aligned}
      \sum_{t = 0}^{t_1} N_t
      &=
      N \sum_{t = 0}^{t_1} 2^t
      =
      N (2^{t_1 + 1} - 1)
      =
      N_{t_1 + 1} - N
      \\
      &\leq
      N_{t_1 + 1}
      =
      2 N_{t_1}
      \leq
      4 [\tilde{K}_0^L + \tilde{K}_1^L + K_1^L(\epsilon)].
    \end{aligned}
  \]
\end{proof}

  \newpage
  \appendix
  \section{Bounds on $\omega$ and its inverse}

\begin{lemma}
  \label{St:Omega-Bounds}
  Let
  $\Map{\omega}{\ClosedOpenInterval{0}{+\infty}}{\RealField}$ and
  $\Map{\omega_*}{\ClosedOpenInterval{0}{1}}{\RealField}$ be the functions
  \[
    \omega(t) \DefinedEqual t - \ln(1 + t),
    \qquad
    \omega_*(\tau) \DefinedEqual -\tau - \ln(1 - \tau).
  \]
  Then, $\omega$ and $\omega_*$ are strictly increasing convex functions on
  their domains.
  For all $t \geq 0$ and all $\tau \in \ClosedOpenInterval{0}{1}$, we have
  \begin{gather}
    \label{Omega-Bounds}
    \frac{t^2}{2 (1 + \frac{2}{3} t)}
    \leq
    \omega(t)
    \leq
    \frac{t^2}{2 + t},
    \\
    \label{Omega-star-Bounds}
    \frac{\tau^2}{2 (1 - \frac{1}{3} \tau)^2}
    \leq
    \omega_*(\tau)
    \leq
    \frac{\tau^2}{2 (1 - \tau)}.
  \end{gather}
  Also, for all $u \geq 0$, the inverse functions can be estimated as follows:
  \begin{gather}
    \label{Inverse-omega}
    \sqrt{2 u} + \tfrac{2}{3} u
    \leq
    \omega^{-1}(u)
    \leq
    \sqrt{2 u} + u,
    \\
    \label{Inverse-omega-star}
    1 - \exp(-\sqrt{2 u} - \tfrac{1}{3} u)
    \leq
    \omega_*^{-1}(u)
    \leq
    1 - \exp(-\sqrt{2 u} - u).
  \end{gather}
\end{lemma}

\begin{proof}
  \ProofPart

  The inequalities~\eqref{Omega-Bounds} and the upper bound in
  \cref{Omega-star-Bounds} were proved in Lemma~5.1.5
  in~\cite{Nesterov-LecturesConvex-18}.
  Let us prove the lower bound in \cref{Omega-star-Bounds}.
  Denote
  \[
    \psi(\tau)
    \DefinedEqual
    \frac{\tau^2}{2 (1 - \frac{1}{3} \tau)^2}.
  \]
  We need to show that $\psi(\tau) \leq \omega_*(\tau)$
  for all $\tau \in \ClosedOpenInterval{0}{1}$.
  Since $\psi(0) = \omega_*(0) = 0$, it suffices to prove the corresponding
  inequality for the derivatives: for all $\tau \in \OpenOpenInterval{0}{1}$,
  \[
    \psi'(\tau) \leq \omega_*'(\tau).
  \]
  But this is easy:
  \[
    \psi'(\tau)
    =
    \frac{\tau}{(1 - \frac{1}{3} \tau)^2}
    +
    \frac{\tau^2}{3 (1 - \frac{1}{3} \tau)^3}
    =
    \frac{\tau}{(1 - \frac{1}{3} \tau)^3}
    \leq
    \frac{\tau}{1 - \tau}
    =
    \omega_*'(\tau).
  \]

  \ProofPart

  Now, let us prove \cref{Inverse-omega}.
  For this, we need to prove that, for all $t \geq 0$,
  \[
    \chi_{2 / 3}(t) \leq t \leq \chi_1(t),
  \]
  where, for any $a \in \OpenClosedInterval{0}{1}$, we denote
  \[
    \chi_a(t)
    \DefinedEqual
    \sqrt{2 \omega(t)} + a \omega(t).
  \]
  Since $\xi_a(0) = 0$, it suffices to prove the inequality for derivatives:
  for all $t > 0$,
  \[
    \xi_{2 / 3}'(t) \leq 1 \leq \chi_1'(t).
  \]
  Using that $\omega'(t) = t / (1 + t)$, we obtain
  \[
    \begin{aligned}
      \xi_a'(t) - 1
      &=
      \frac{\omega'(t)}{\sqrt{2 \omega(t)}} + a \omega'(t) - 1
      =
      \frac{t}{1 + t}
      \Bigl( \frac{1}{\sqrt{2 \omega(t)}} + a \Bigr) - 1
      \\
      &=
      \frac{t}{1 + t}
      \Bigl( \frac{1}{\sqrt{2 \omega(t)}} - \frac{1 + (1 - a) t}{t} \Bigr).
    \end{aligned}
  \]
  Thus, we need to check if
  \[
    \frac{1}{t}
    \leq
    \frac{1}{\sqrt{2 \omega(t)}}
    \leq
    \frac{1 + \frac{1}{3} t}{t},
  \]
  or, equivalently, if
  \[
    \frac{t^2}{2 (1 + \frac{1}{3} t)^2}
    \leq
    \omega(t)
    \leq
    \frac{t^2}{2}.
  \]
  But these inequalities easily follow from \cref{Omega-Bounds}.

  \ProofPart

  Finally, let us prove \cref{Inverse-omega-star}.
  These inequalities are equivalent to, for all $u \geq 0$,
  \[
    \sqrt{2 u} + \tfrac{1}{3} u
    \leq
    -\ln\bigl( 1 - \omega_*^{-1}(u) \bigr)
    \leq
    \sqrt{2 u} + u.
  \]
  Denote
  \[
    \chi_a(\tau)
    \DefinedEqual
    \sqrt{2 \omega_*(\tau)} + a \omega_*(\tau).
  \]
  Then, we need to prove that, for all $\tau \in \ClosedOpenInterval{0}{1}$, we
  have
  \[
    \chi_{1 / 3}(\tau)
    \leq
    -\ln(1 - \tau)
    \leq
    \chi_1(\tau),
  \]
  or, equivalently, in view of the definition of $\omega_*$,
  \[
    \chi_{1 / 3}(\tau)
    \leq
    \tau + \omega_*(\tau)
    \leq
    \chi_1(\tau).
  \]
  This is the same as
  \[
    \chi_{-2 / 3}(\tau) \leq \tau \leq \chi_0(\tau).
  \]
  Since $\chi_a(0) = 0$, it suffices to prove these inequalities for the
  derivatives: for all $\tau \in \OpenOpenInterval{0}{1}$,
  \[
    \chi_{-2 / 3}'(\tau) \leq 1 \leq \chi_0'(\tau),
  \]
  Since $\omega_*'(\tau) = \frac{\tau}{1 - \tau}$, we obtain
  \[
    \begin{aligned}
      \chi_a'(\tau) - 1
      &=
      \frac{\omega_*'(\tau)}{\sqrt{2 \omega_*(\tau)}}
      +
      a \omega_*'(\tau)
      -
      1
      =
      \frac{\tau}{1 - \tau}
      \Bigl( \frac{1}{\sqrt{2 \omega_*(\tau)}} + a \Bigr)
      -
      1
      \\
      &=
      \frac{\tau}{1 - \tau}
      \Bigl(
        \frac{1}{\sqrt{2 \omega_*(\tau)}}
        -
        \frac{1 - (1 + a) \tau}{\tau}
      \Bigr).
    \end{aligned}
  \]
  Thus, we need to check if
  \[
    \frac{1 - \tau}{\tau}
    \leq
    \frac{1}{\sqrt{2 \omega_*(\tau)}}
    \leq
    \frac{1 - \frac{1}{3} \tau}{\tau},
  \]
  or, equivalently, if
  \[
    \frac{\tau^2}{2 (1 - \frac{1}{3} \tau)^2}
    \leq
    \omega_*(\tau)
    \leq
    \frac{\tau^2}{2 (1 - \tau)^2}.
  \]
  But this readily follows from \cref{Omega-star-Bounds}.
\end{proof}

  \printbibliography
\end{document}